\documentclass[reqno]{amsart}
\usepackage{hyperref}
\usepackage{graphicx}

\usepackage{amsmath}
\usepackage{amssymb}
\usepackage{graphicx}
\usepackage{xcolor}
\usepackage{comment}
\usepackage{makecell}

\AtBeginDocument{{\noindent\small
\emph{Manuscript}.\newline 
}
\vspace{8mm}}

\begin{document}
\title[Flow in porous medium layered over inclined impermeable bed]{Modeling of groundwater flow in
porous medium layered over inclined impermeable bed}
%\title[Flow in porous medium layered over inclined impermeable bed]
%{Mathematical Model 
%of groundwater flow in porous medium layered over inclined impermeable bed and some analytical tools of its study}

\author[P. Girg, L. Kotrla\hfil \dots \hfilneg]
{Petr Girg,  Luk\'{a}\v{s} Kotrla}

\address{Petr Girg \newline
Department of Mathematics and NTIS,
Faculty of Applied Scences, University of West Bohemia,
Univerzitn\'{\i} 8, CZ-301\,00~Plze\v{n}, Czech Republic}
\email{pgirg@kma.zcu.cz}

\address{Luk\'a\v{s} Kotrla \newline
Department of Mathematics and NTIS,
Faculty of Applied Scences, University of West Bohemia,
Univerzitn\'{\i} 8, CZ-301\,00~Plze\v{n}, Czech Republic}
\email{kotrla@ntis.zcu.cz}

\dedicatory{In memory of Professor John~W.~Neuberger with admiration}

%\thanks{Submitted \dots. Published \dots, \dots.}
\keywords{Porous medium; filtration; nonlinear Darcy's law;
\hfill\break\indent pressure-to-velocity power law; $p$-Laplacian}

\begin{abstract}
 We propose a new mathematical model of groundwater flow in porous medium layered over inclined impermeable bed.
 In its full generality, this is a free-surface problem. To obtain analytically tractable model, we use generalized Dupuit-Forchheimer assumption for inclined impermeable bed.
 In this way, we arrive at parabolic  partial differential equation which is a generalization of the classical Boussinesq equation.  
 Novelty of our approach consists in considering nonlinear constitutive law of the power type.
 Thus introducing $p$-Laplacian-like differential operator into the Boussinesq equation. Unlike in the classical case of the Boussinesq equation, the convective term cannot be set aside from the main part of the diffusive term and remains incorporated within it. 
 
In the sequel of the paper, we analyze qualitative properties of the stationary solutions of our model. In particular, we 
study existence and regularity of weak solutions for the following boundary value problem
\begin{equation*}
\begin{aligned}
	&
	- 
	\frac{\rm d}{{\rm d} x}
	\left[
	(u(x) + H) \left|\frac{{\rm d} u}{{\rm d} x}(x) \cos(\varphi) + \sin(\varphi) \right|^{p - 2}
	\left(\frac{{\rm d} u}{{\rm d} x}(x) \cos(\varphi) + \sin(\varphi)\right)
	\right]
	\\
	&
	\begin{aligned}
	&
	= 
	f(x)\,, & \qquad\qquad x \in (-1,1)\,,\\ 
	&
	u(-1) = u(1) = 0\,,&
	\end{aligned}
\end{aligned}
\end{equation*}
where $p>1$, $H>0$, $\varphi\in (0, \pi/2)$, $f\geq 0$,  $f\in L^{1}(-1,1)$.
In the case of $p>2$,
we  study validity of Weak and Strong Maximum Principles as well.
We use methods based on the linearization of the 
$p$-Laplacian-type problems in the vicinity of known solution,
error estimates,
and analysis of Green's function of the linearized problem. 
\end{abstract}

\maketitle
\numberwithin{equation}{section}
\newtheorem{theorem}{Theorem}[section]
\newtheorem{remark}[theorem]{Remark}
\newtheorem{definition}[theorem]{Definition}
\newtheorem{lemma}[theorem]{Lemma}
\newtheorem{corollary}[theorem]{Corollary}
\newtheorem{proposition}[theorem]{Proposition}
\newtheorem{counterexample}[theorem]{Counterexample}
\newtheorem{openquestion}[theorem]{Open Question}
\allowdisplaybreaks

\section{Introduction} \label{s:Intro}
Hand in hand with global warming 
(whether man made or not), global water cycle intensifies and hydrological extremes (such as heavy precipitation events and local floods) may occur more frequently, see, e.g., \cite{hess-18-2735-2014, HUNTINGTON200683, Trenberth2011123, Westra2014522}. Thus, further research and development of more effective drainage systems is needed. A~typical situation which 
frequently appears in this context is water flow in a porous medium layered over a slopping impermeable bed. In practice, this situation can be encountered, e.g., in highway and railway drainage~\cite{GhataoraRushton2012,YoungsRushton2009}, buried streams through coarse porous
media and valley fills \cite{HosseiniJoy2007, Sedghi-Asl2016, Sedghi-Asl20141D}, 
water seepage through soil between parallel ditches in irrigated/drained sloping lands \cite{Chapman1980, 
Childs71, 
SchmidLuthin1964,
SinghRaiRamana1991, WoodingChapman1966} and/or groundwater flow in inclined phreatic aquifers \cite{Bansal2015, BansalDas2011, BarlettPorporato2018, YOUNGS1990201}. In its full generality, this flow configuration leads to a free boundary value problem, which is often very difficult to analyze and/or use in simulations, see, e.g., \cite{BaiocchiComincioliMagenesPozzi1973}. 
For easier application in practice, several simplified mathematical models  were proposed and studied in a number of papers, see, e.g., 
\cite{BarlettPorporato2018, Chapman1980,
Childs71,
DalyPorporato2004,
Loaiciga2005, 
WoodingChapman1966}
and references therein. Most of these models study groundwater flow in a soil layered over impermeable bed and use linear constitutive law, known as Darcy's law (see \cite{Darcy1856}, and also, e.g.,
\cite{AravinNumerov,
Bear1972, BenediktGirgKotrlaTakac2018, Harr2012groundwater, Polubarinova-Kochina1952, Scheidegger1960} for further discussions). However, it has been established by numerous complex laboratory experiments and by many in-field observations that the Darcy's law is not satisfactory for certain materials and flow regimes, see, e.g.,
\cite{Forchheimer1901, Izbash1931, Kroeber1884, Missbach1936, Smreker1878, Smreker1879, SoniIslamBasak1978}.
However, typical materials used in modern drainage systems comprise of coarse porous media such as gravel or geosynthetic materials and it turns out that movement of water in such materials is more accurately described by nonlinear constitutive laws such as power type law, sometimes called Smreker-Izbash-Missbach law, or polynomial type law, known as Forchheimer's law. 
Empirical studies on these materials can be found in
\cite{Sedghi-Asl2014}
(gravels),
\cite{BORDIER2000174} (geosynthetic materials),
\cite{Eck2012}
(porous asphalt), 
\cite{Koohmishi2019, SchmidtEtAl2017} (railway ballast material).
Moreover, it has been experimentally established that 
the groundwater obeys nonlinear constitutive law of power type
in low permeability materials such as certain  sandstones, fine sands, clays or certain types of soils, see, e.g., 
\cite{King1898, SoniIslamBasak1978, Zunker1920}.

The purpose of this paper is twofold. At first we propose an improved mathematical model of  water flow in a~porous medium obeying power law, layered over a slopping impermeable bed in Section~\ref{sec:model}. Then we present some mathematical tools suitable for its study in Section~\ref{sec:prop:w:sol}. These tools are mainly based on linearization of quasilinear operators of the $p$-Laplacian type. In particular, we prove validity of Strong Maximum Principle in certain situations that can be met in analyzing real world situations. 
Concluding remarks with
a table summarizing results of the papers are in Section~\ref{sec:concrem}.
Finally, technical parts of some proof from Section~\ref{sec:prop:w:sol} can be found in the Appendix~\ref{app:a}.

We use usual notation throughout the paper such as, e.g.,  
$L^p(-1,1)$, $p\geq 1$, stands for spaces 
of Lebesgue integrable functions,
$C[a,b]$, $C(a,b)$ 
with $a<b$ stand for spaces of continuous functions on respective intervals $[a,b]$ and $(a,b)$, $AC[-1,1]$ stands for the space of absolutely continuous functions on $[-1,1]$. Analogously, $C^k[a,b]$ and $C^k(a,b)$, $k\in\mathbb{N}$, stands for spaces of $k$-times continuously differentiable functions, $C^{\infty}(a,b)=\bigcap_{k\in\mathbb{N}} C^k(a,b)$ and 
$C_{\mathrm{c}}^{\infty}(-1,1)$ stands for the space of smooth functions $u$ with compact support $\mathop{\mathrm{supp}}u\subset (-1,1)$.
By 
$W_0^{1,p}(-1,1)$, $p>1$, we denote the space of functions $u$ from $L^p(-1,1)$ having its distributional derivatives in $L^p(-1,1)$ and satisfying 
$u(-1)=u(1)=0$ in the sense of traces. Let us note that the boundary conditions are satisfied in the classical sense, since any $u\in W_0^{1,p}(-1,1)$
has a representative in $ AC[-1,1]$ (see, e.g., \cite[Thm.~2.1.4]{Ziemer89}).  
We denote the norm of a Banach space $X$ by $\|\cdot\|_X$,
with the only exception $L^{\infty}(-1,1)$ whose corresponding norm is denoted by $\|\cdot\|_{\infty}$ for brevity.
Throughout the paper, we will also use the following standard notation 
$$
u^+\stackrel{\mathrm{def}}{=}\max\{u,0\}\,,\quad
u^-\stackrel{\mathrm{def}}{=}\max\{-u,0\}\,,
$$
for any real function
$u$ (recalling that
$u = u^+ - u^-$).

\section{Mathematical model of water flow in porous medium layered over inclined impermeable bed}
\label{sec:model}
\subsection{Initial physical considerations}
Generally, water tends to flow from places with higher potential energy to places with lower potential energy.
%In general, groundwater flows from higher elevated places to lower elevated places (if capillarity effects are neglected). 
While the water is flowing, its potential energy is transformed into kinetic energy, which is further dissipated by viscosity forces and friction with the porous medium. 
Since the actual velocity of the groundwater highly oscillates in the channels in the porous medium, one must consider bulk motion of the water within a sufficiently large control volume in the porous medium. This averaged velocity $\vec{v}$ is defined by means of specific discharge $\vec{q}$
by formula $\vec{v}=\vec{q}/n$, where
the specific discharge (vector) $\vec{q}$ takes the direction of the flow and its magnitude is defined as the volume of water flowing per unit time through a unit cross-sectional area normal to the direction of flow, and $n$ is effective porosity of the medium
see, e.g., 
\cite[p.~121]{Bear1972} for detailed explanation. On this macroscopic level, the process of transformation and dissipation of energy can be described in the following way. The total mechanical energy per volume in a control volume of water is the sum of gravitational potential energy, pressure energy, and kinetic energy
$$
E_{T}=z \varrho g+P+\frac{1}{2} \varrho v^{2}\,,
$$
where $v$ stands for the magnitude of averaged velocity of the flow in the control volume, $\varrho$ is the water density, 
$P$ pressure, 
$z$ elevation of the control volume from the datum, 
$g$ gravitational acceleration, see, e.g., \cite{ZekaiSen1995}. For incompressible liquid such as water, one can equivalently consider another quantity called total head
$$
h_{T}\stackrel{\mathrm{def}}{=}\frac{E_T}{\varrho g}=z+\frac{P}{\varrho g}+\frac{1}{2 g} v^{2}
$$
which can be directly measured in practice, e.g., by using observation wells or by so called piezometers, see, e.g., \cite[p.~63]{Bear1972} or \cite[pp.~129--132]{ZekaiSen1995}. As it was already mentioned, groundwater is loosing its total energy (or equivalently total head) while flowing due to viscous forces and friction with porous medium. Thus, its total energy decreases in the direction of the flow. 
%is maximally of the order of a meter per %day (that is of the order $0.00001 %\mathrm{~m} / \mathrm{s}$ ), 
The term corresponding to kinetic energy is negligible and can be dropped for the average velocities of the groundwater flow encountered in real situations, see, e.g., 
\cite[p.~5]{Harr2012groundwater} or
\cite[pp.~40--43]{ZekaiSen1995}. 
In this way, we obtain piezometric head
$$
h\stackrel{\mathrm{def}}{=}z+\frac{P}{\varrho g}
$$
which is the state variable in the mathematical models of the water flow in the underground. On the other hand, the specific discharge is the flux quantity.
The constitutive law relating this two quantities, quantitatively describes the rate of dissipation of the energy along the flow path. 

\subsection{Constitutive law}
In practice, the constitutive law is obtained empirically from experimental data for given porous medium and fluid.  
Linear Darcy's law relates groundwater flux to the piezometric head loss per length according to the following formula
$$
v = c\, \frac{\triangle h}{\triangle L}\,,
$$
where $c>0$ is a constant to be determined from measured data,  
$\triangle h$ is the difference of the piezometric head measured at two distinct locations distance $\triangle L$ apart. 
This formula was established experimentally for filtration of water through sand by Henry Darcy \cite{Darcy1856} in 1856 and soon it became widely popular in mathematical models of groundwater flow due to its simplicity (and still reasonable accuracy).
Unfortunately, it was later found that it has limited range of its validity in coarse grained media (such as gravels), see, e.g., \cite{Forchheimer1901, Izbash1931, Kroeber1884, Missbach1936, Smreker1878, Smreker1879} as well as in media with very low permeability (such as clays, certain soils, sandstones), see, e.g., \cite{King1898, Zunker1920}. For thorough surveys and discussions of this and other constitutive laws and various criteria of their validity, see, e.g., \cite{AravinNumerov, 
Bear1972, Scheidegger1960,
ZekaiSen1995,
SoniIslamBasak1978}. It follows from these discussions and the above mentioned papers that the power type law
\begin{equation}
\label{power:law}
v = c\, \left(\frac{\triangle h}{\triangle L}\right)^m
\end{equation}
with  constants $c, m >0$ to be empirically determined,
turned out to be simple but flexible enough to fit with most experimental data obtained for various porous media. For $m=1$, the power law coincides with the Darcy's law. 

Note that the
constitutive laws are inferred from experiments made on one dimensional flow and that the averaged velocity is (practically) constant within the sample of material and during the time of the measurement. 
However, groundwater flow in the real world is three dimensional, in general,  and the physical quantities $\vec{v}$ and $h$ are usually functions of spatial variables and time. In case
of the homogeneous and isotropic porous medium, the three-dimensional constitutive law can be inferred from the one-dimensional one in a straightforward manner, taking into account that
the averaged velocity takes the opposite direction of the gradient of the piezometric head and no 
flow occurs if the gradient of the piezometric head is zero. In this way, we obtain
\begin{equation}
\label{vec:power:law}
\vec{v} = 
\begin{cases} \vec{0} & \text { for } \nabla h=\vec{0}\,, \\ 
- c\, \left|\nabla h\right|^{m-1} \nabla h  & 
\text { for } 
\nabla h\neq\vec{0}\,,
\end{cases}
\end{equation}
where $\nabla h$ stands for the spatial gradient of the piezometric head, $c, m>0$ are constants as in \eqref{power:law}.

Note that the power law \eqref{power:law} is not the only type of nonlinear laws used in practice. For thorough surveys of other important types of nonlinear constitutive laws, see, e.g., \cite{AravinNumerov, Bear1972, Scheidegger1960, ZekaiSen1995}. Furthermore, a discussion of development of constitutive laws for fluid flows in porous media from the perspective of history of science can be found, e.g., in  \cite{AravinNumerov, BenediktGirgKotrlaTakac2018}. 
%Measurement power-law railway ballast material:\cite{Koohmishi2019}.

\subsection{Free surface problem and Dupuit-Forchheimer assumption}

\phantom{b}\newline The groundwater flow
with a free-surface upper boundary in a porous medium layered over impermeable  bed is  very challenging problem, since a part of the boundary of the domain is not known 
a~priori and thus it is one of the unknowns. Rigorous formulation and mathematical treatment of some important cases 
of problems with  free-surface upper boundary can be found, e.g., in~\cite{BaiocchiComincioliMagenesPozzi1973}. Since these types of problems are frequently encountered in engineering, it was as early as in the middle of the nineteenth century when French engineer J.~Dupuit~\cite{Dupuit1863} (cf.~also \cite{Dupuit1848} for similar approach to flow in open channels)
published a relatively simple method how to find approximate solutions of these type of problems and used this method to study groundwater flow with free surface towards a fully penetrating well.
His method was based on observation that the maximum slope of the upper free surface is very small, (typically $\triangle h/ \triangle L$ is of order $0.001$). This lead him to the following two simplifying assumptions:
\begin{itemize}
\item[{(D1)}] 
the groundwater flow is horizontal (piezometric head is constant in vertical direction),
\item[{(D2)}] 
the groundwater flow is proportional to the gradient of the piezometric head.
\end{itemize}

\begin{figure}
\begin{picture}(350,290)
    \put(0,0){\includegraphics[width=0.95\textwidth]{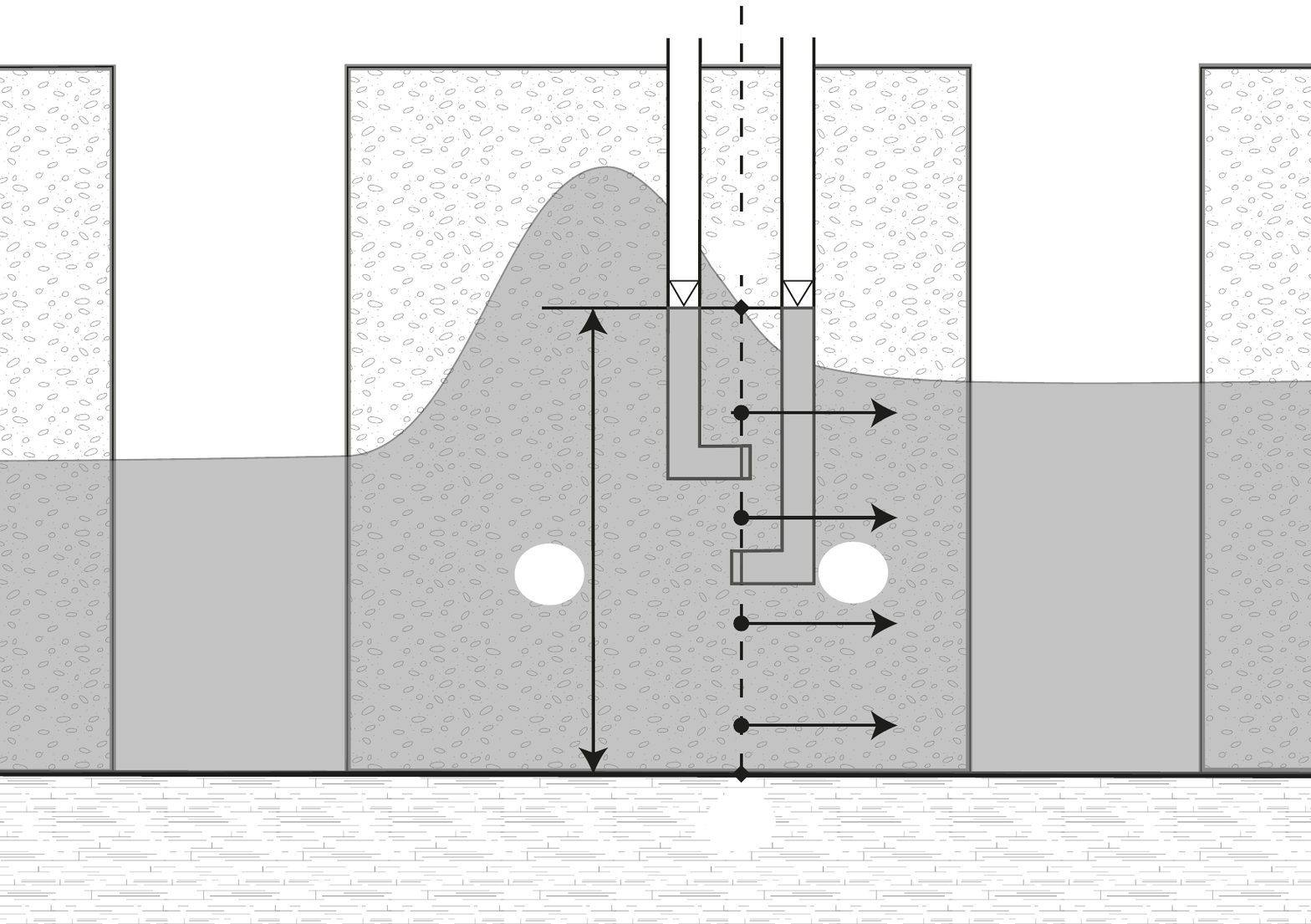}}
    \put(189,24){$\mathrm{A}$}
    \put(189,172){$\mathrm{B}$}
    \put(220,89){$\vec{v}$}
    \put(140,88){$h$}
    \put(140,250){\vector(0,-1){20}}
    \put(150,250){\vector(0,-1){20}}
    \put(160,250){\vector(0,-1){20}}
    \put(170,250){\vector(0,-1){20}}
    \put(145,260){rain}
\end{picture}
\caption{Classical Dupuit-Forchheimer assumption (D1). The flow is horizontal with the velocity distribution profile and the piezometric head being constant along vertical line segment $\mathrm{A}\mathrm{B}$. }
\label{fig:head}
\end{figure}

Since the capillarity effects are also neglected in his approach, 
the hydrostatic pressure must be equal to atmospheric pressure at the free-surface.
Therefore, he obtained the relation 
$$h(x, y, z, t) = z\,,
$$
for any point $(x,y,z)$ at the free surface at any time $t$, see Fig.~\ref{fig:head}.
It then follows from the assumption (D1) 
that 
$
h(x,y,z,t)=h(x,y,t)
$ (with a little abuse of notation) and
the height $z$ of the free-surface 
boundary above the impermeable layer is then $z=h(x,y,t)$.
Thus it reduces the original 3D problem with unknown free-surface upper boundary to a problem in the $xy$-plane, where the surveyed area is known in advance, see, e.g., \cite{AravinNumerov,
Bear1972, Harr2012groundwater}
for derivation of the corresponding equations in modern notation. 

Independently, similar research was performed by Austrian engineer
P.~Forch\-heimer~\cite{Forchheimer1886}, where 
groundwater flows towards a system of fully penetrating wells, or 
towards a fully penetrating slot of finite length were studied. The approach by Dupuit~\cite{Dupuit1863} and Forchheimer~\cite{Forchheimer1901} turned~out to be very useful and enabled to solve many practical problems arising from engineering. However, one has to bear in mind that it is only an approximation of the original problems.
There is a thorough discussion in \cite{Bear1972, Harr2012groundwater} of various cases, for which Dupuit-Forchheimer assumption leads to a reasonable approximation of the free surface problem as well as the cases for which it yields a poor approximation. Results obtained using the Dupuit-Forchheimer assumption were compared to experimentally measured data, e.g., in \cite{YOUNGS1990201}.

In \cite{Forchheimer1901, Zhukovskii1889}, these ideas were extended to groundwater flow in natural coarse grained media such as gravel, where nonlinear constitutive laws better fits experimental data. 
In both papers \cite{Forchheimer1901, Zhukovskii1889}, authors used assumption (D1) with 
nonlinear power law \eqref{power:law} (there, written in another equivalent forms) to solve a problem of groundwater flow towards a fully penetrating well, and also towards an infinite fully penetrating ditch in case of \cite{Forchheimer1901}. 

\subsection{Dupuit-Forchheimer assumption on the sloping porous medium}
\label{par:Dupuit:Forchheimer:sloping}
The situation becomes even more complicated when the impermeable 
bed (upon which the porous medium is resting) is inclined. There are two main approaches how to adopt Dupuit-Forchheimer assumption to this situation. The simpler approach is based on (D1),
that is, the flow is horizontal despite the impermeable bed is inclined. Thus it turns out that this approach is reasonable only for very small inclinations.
It was already in 1848, when
Dupuit~\cite{Dupuit1848} 
introduced this approach to study 
water flow in an inclined open channel and later adopted in 
Boussinesq~\cite{Boussinesq1904} for the groundwater flow in a porous medium over an inclined impermeable bed. 
For its simplicity, it is widely used in solving engineering problems
see, e.g., monographs
\cite{Bear1972, Harr2012groundwater} and 
also
\cite{Bansal2015,
BansalDas2009, BansalDas2010, Chapman1980,
SchmidLuthin1964}. 

The second approach suggested already by
Boussinesq~\cite[pp.~252--260]{Boussinesq1877} in 1877, and later promoted by Childs~\cite{Childs71}, 
is based on the so-called extended Dupuit–Forch--heimer assumption,
that is,
\begin{itemize}
\item[{(eD1)}] the groundwater flow is parallel to the inclined impermeable bed and the piezometric head is constant in the normal direction to the inclined impermeable bed.
\end{itemize}  
The latter approach was used, e.g., 
in~\cite{BansalDas2011, Towner1975, VerhoestTroch2000}.

Very interesting discussions about the validity, use, and comparison of these two approaches can be found, e.g., in \cite{WoodingChapman1966} and introductions of papers \cite{BansalDas2011, BaruaMazumdar2020}. 
In Towner~\cite{Towner1975}, the second approach (flow parallel to inclined impermeable bed) was used to study groundwater flow between parallel ditches in case of uniform rainfall. It was found that the calculated water table heights are in much better agreement with experimental data than those published in \cite{SchmidLuthin1964} and calculated using the first approach assuming horizontal flow. {\it Based on this, we adopt the second approach in the model developed in this paper}.

\subsection
{Physical assumptions of our model}
\begin{figure}
\begin{picture}(350,290)
 \put(0,0){\includegraphics[width=\textwidth]{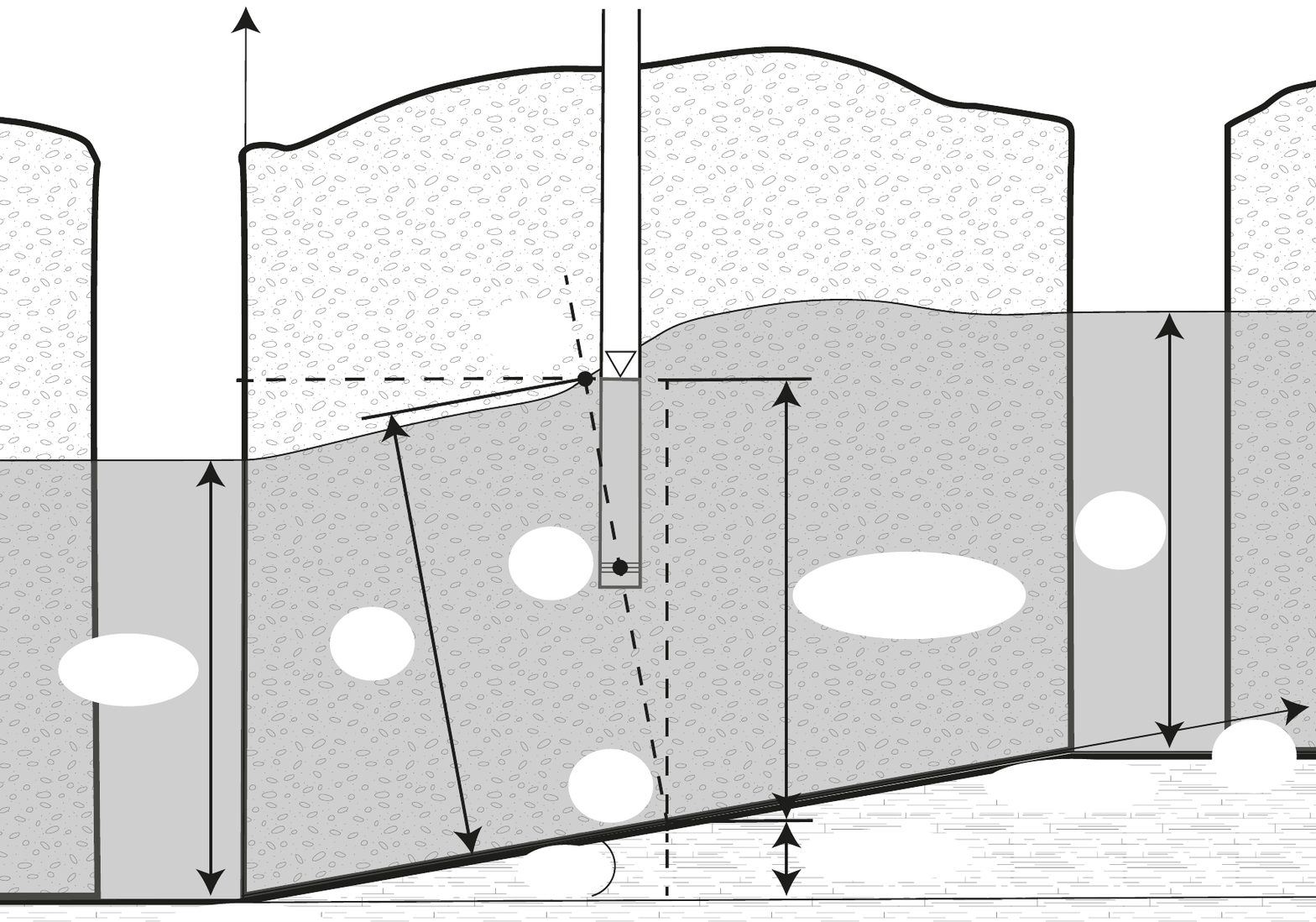}}
\put(151,13){$\varphi$}
\put(20,66){$h(-1)$}
\put(297, 105){$h(1)$}
\put(99, 72){$\widehat{h}$}
\put(53, 145){$h$}
\put(232, 87){$\widehat{h}\,\cos\varphi$}
\put(225, 16){$x\sin\varphi$}
\put(340, 42){$x$}
\put(55, 233){$z$}
\put(140, 158){$A$}
\put(163, 34){$B$}
\put(148, 94){$C$}
\put(50, -5){$x=-1$}
\put(280, 35){$x=1$}
\put(180,260){\vector(0,-1){20}}
\put(190,260){\vector(0,-1){20}}
\put(200,260){\vector(0,-1){20}}
\put(210,260){\vector(0,-1){20}}
\put(185,270){rain}
\end{picture}
\caption{Geometric relations between $h$, $\widehat{h}$ and $\varphi$. Due to (eD1) Piezometric head $h$ is assumed to be constant along the line segment $AB$. The fictive piezometer is measuring piezometric head at point $C$.}
\label{fig1}
\end{figure}

We assume that there is an inclined layer of impermeable bed (such as, e.g., bedrock or impermeable geosynthetic material) covered by a parallel layer of permeable material (such as, e.g., permeable rock, soil, sand, gravel or some geosynthetic porous material). 
The groundwater moves
within the permeable layer in the saturated zone only
(see, e.g., \cite[Sec.~1.1.2, p.~2 ]{Bear1972}). We assume that the groundwater has a~free surface, capillarity effects above the free surface are neglected, and that the~extended~Dupuit-Forchheimer assumption (eD1) for the free surface is valid, i.e.,  the piezometric head $h(\,\cdot\,,t)$ is constant in the direction perpendicular to the impermeable bed at any given time $t$. For simplicity, we further assume that there is no bulk flow in the $y$ direction.
Thus the piezometric head is a function of $x$ and $t$ only and given $(x,t)$, the free surface is at the distance $h(x,t)$ from the 
impermeable bed, see Fig.~\ref{fig1}. 
We assume that the groundwater flow through porous medium forming the permeable layer
is governed by 
the power law \eqref{vec:power:law} which in our situation is reduced to
\begin{align*}
	v(x,t)
	& =
	-c 
	\left|
	\frac{\partial h}{\partial x} (x,t) 
	\right|^{m - 1}
	\frac{\partial h}{\partial x} (x,t)\,. 	
\end{align*}

\begin{comment}
\begin{figure}
\centering
    \includegraphics[width=0.95\textwidth]{Naklonena_hladina}
\label{fig2}    
\end{figure}
\end{comment}

\begin{figure}
\centering
    \includegraphics[width=0.95\textwidth]{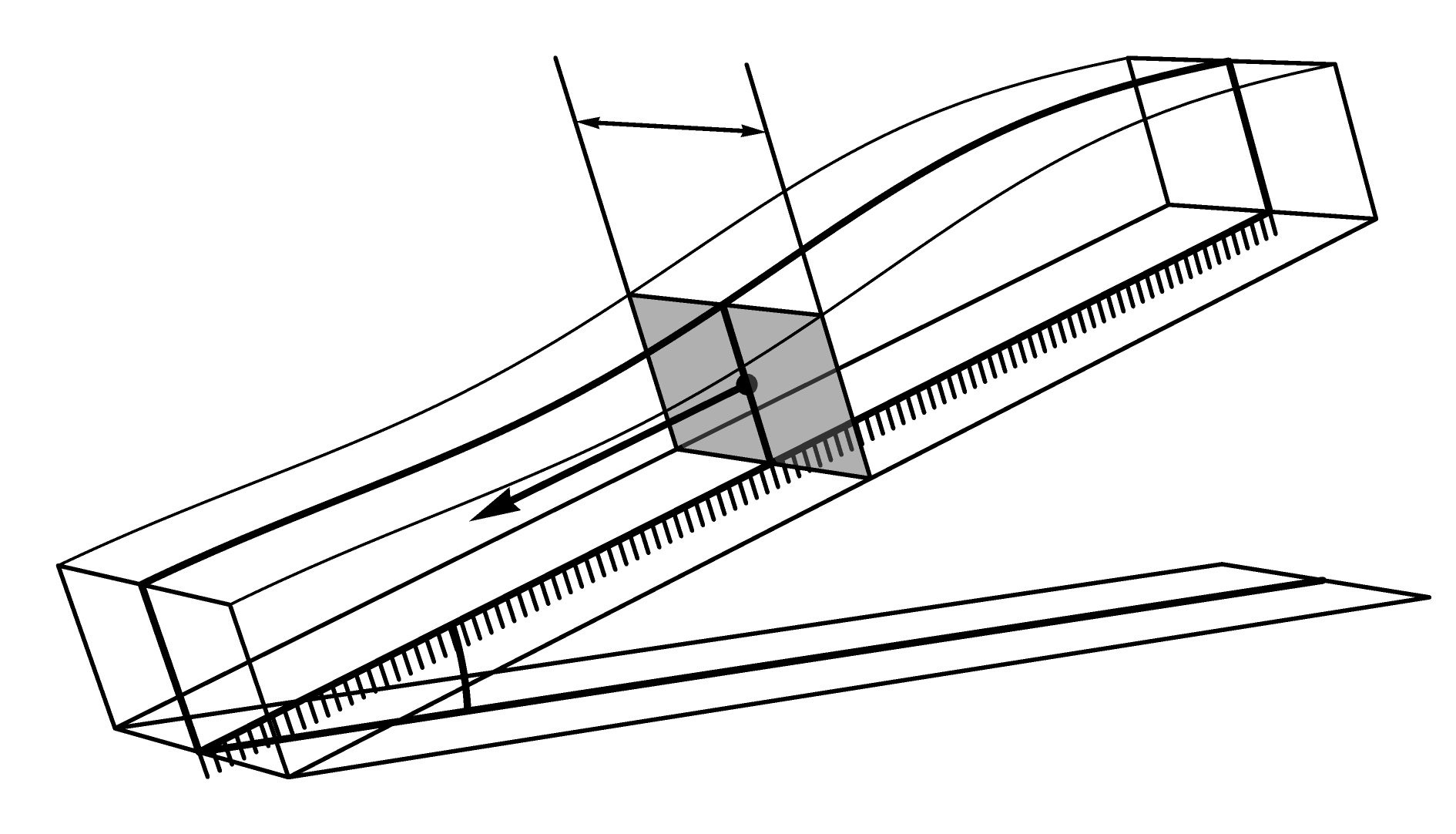}
\put(-195,173){unit}
\put(-195,156){width}
\put(-145,160){
\fcolorbox{white}{white}{free surface}}
\put(-130, 120){
\fcolorbox{white}{white}{impermeable layer}}

\centering 
(a) the 3D geometric configuration.

\centering
    \includegraphics[width=0.95\textwidth]{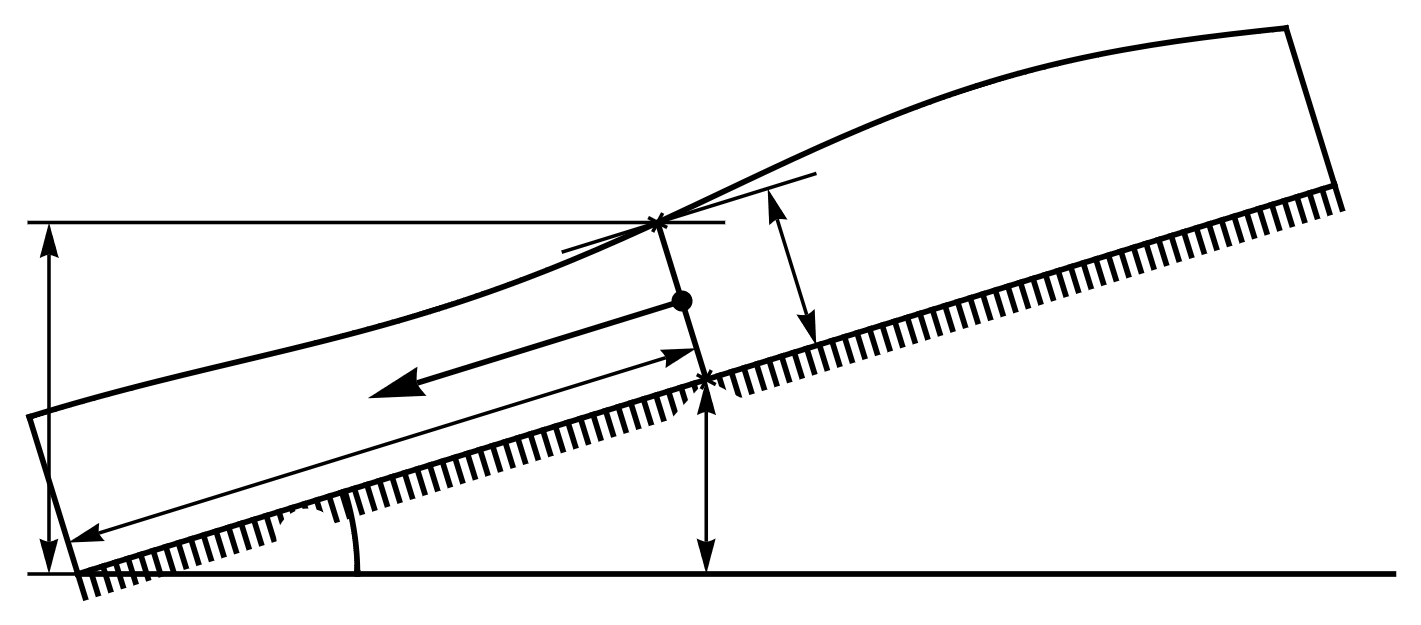}
\put(-271,18){$\varphi$}
\put(-165,30){$x \sin\varphi $}
\put(-145,85){$\widehat{h}(x,t)$}
\put(-268,45){$x$}
\put(-220,72){$v$}
\put(-190,100){$A$}
\put(-170,65){$B$}
\put(-325,70){$h(x,t)$}
\put(-145,125){\fcolorbox{white}{white}{free surface}}
\put(-90,80){\fcolorbox{white}{white}{impermeable}}
\put(-80,67){\fcolorbox{white}{white}{layer}}

\centering (b) 2D longitudinal-vertical section with state and flow variables marked.

\caption{State and flow variables in the case of the flow in porous medium layered over inclined impermeable layer under generalized Dupuit-Forchheimer assumption (eD1). 
The averaged velocity is assumed to be constant on the line segment AB and parallel to the inclined bedrock. Here,  
$h(x,t)= \widehat{h}(x,t)\cos\varphi + x\sin\varphi$,
$v=v(x,t)=|\partial h/\partial x|^{m-1} \partial h/\partial x$,
and the groundwater flow per unit width is $Q_{\mbox{\tiny p.u.w.}}(x,t)=\widehat{h}(x,t)v(x,t)$.}
\label{fig3}    
\end{figure}

Here, $v(x,t)$ is flow velocity, and $h(x,t)$ is piezometric head
measured from the chosen horizontal reference level, see Fig.~\ref{fig1}.
Then, the flux 
per unit width in the saturated zone 
(denoted by 
$Q_{\mbox{\tiny p.u.w.}}$)
at any position $x$
(see Fig.~\ref{fig1}) and any time $t$ 
is
\begin{equation}
\label{PowerLawOurCase}
	Q_{\mbox{\tiny p.u.w.}}(x,t) 
	= 
	\widehat{h}(x,t) v(x,t) 
	= 
	- c\,
	\widehat{h}(x,t)
	\left|
	\frac{\partial h}{\partial x} (x,t) 
	\right|^{m - 1}
	\frac{\partial h}{\partial x} (x,t)\,.
	%\left|\frac{\partial h}{\partial x} \cos(\varphi) + \sin(\varphi) \right|^{p - 2}\left(\frac{\partial h}{\partial x} \cos(\varphi) + \sin(\varphi)\right)\,.
\end{equation}
Here, $\widehat h(x,t)$ is thickness of the saturated zone measured perpendicularly to the inclined impermeable bed, see Fig.~\ref{fig1}.
By the extended Dupuit-Forchheimer assumption (eD1), the piezometric head $h(x,t)$ is constant along planes perpendicular to impermeable bed. Thus 
\begin{equation}
\label{ThicknessToHead}
h(x,t) = 
\widehat{h}(x,t) \cos\varphi +
x \sin\varphi\,,
\end{equation}
where $\varphi$ is angle of inclination, see Fig~\ref{fig1}.

The equation of continuity for our geometric configuration 
reads as follows
\begin{equation}
\label{ContinuityEquation}
	\frac{\partial \widehat{h}}{\partial t}(x,t) + \frac{\partial Q_{\mbox{\tiny p.u.w.}}}{\partial x}(x,t) = r(x,t)\cos \varphi\,,
\end{equation}
since $\widehat{h}(x,t)$ is measured in perpendicular direction to the bedrock and 
thus represents the amount of groundwater in the saturated zone 
per length $x$ and unit width.
The source function $r$ (rain) is given in the vertical direction. Thus $r(x,t)$ has to be multiplied by $\cos\varphi$.
Combining equations
\eqref{PowerLawOurCase}, \eqref{ThicknessToHead} and \eqref{ContinuityEquation}, 
we obtain
\begin{equation}
\label{eq_h_parabolic}
\begin{aligned}
	&\frac{\partial \widehat{h}}{\partial t}(x,t) \\
	& - 
	c \frac{\partial}{\partial x}
	\left[
	\widehat{h}(x,t) \left|\frac{\partial \widehat{h}}{\partial x}(x,t) \cos\varphi + \sin\varphi \right|^{m - 1}
	\left(\frac{\partial \widehat{h}}{\partial x}(x,t) \cos\varphi + \sin\varphi\right)
	\right]\\
	& = r(x,t)\cos\varphi\,.
\end{aligned}
\end{equation}

\subsection{Initial-boundary value problem for a mathematical model of groundwater flow between infinite parallel ditches.}

From now on, we will make use of the parameter 
$p\stackrel{\mathrm{def}}{=} m+1>1$ instead of $m>0$, in order to follow standard notation used by ``$p$-Laplacian community'', whereas $m$ is mostly used in hydrological literature. As usual, $p' > 1$ denotes the conjugate exponent to $p>1$, i.e., $1/p + 1/p' = 1$,  in the sequel of the paper.

Let the constants $H_{-1}, H_{1}\geq 0$ denote the water level measured from the impermeable bed in the left ($x = -1$) and right ($x = 1$) ditch, respectively. 
Let the function $h_0(x)$, $x \in (-1,1)$, describes the initial state of the water level between the ditches.
The water level is measured perpendicularly to the impermeable bed here 
(see Section~\ref{par:Dupuit:Forchheimer:sloping} for more details). To simplify notation, we also define
 $f(x,t)\stackrel{\mathrm{def}}{=} r(x,t) \cos\varphi$.
Then, based on \eqref{eq_h_parabolic}, we obtain initial boundary value problem
\begin{equation}
\label{eq:math:model:final}
\begin{aligned}
	&
	\frac{\partial \widehat{h}}{\partial t}(x,t)
    - 
	c \frac{\partial}{\partial x} 
	\left[
	\widehat{h}(x,t) \left|\frac{\partial \widehat{h}}{\partial x}(x,t) \cos\varphi + \sin\varphi \right|^{p - 2}
	\left(\frac{\partial \widehat{h}}{\partial x}(x,t) \cos\varphi + \sin\varphi\right)
	\right]
	\\
	& 
	\begin{aligned}
	& = f(x,t)\,,& & \qquad\qquad x \in (-1,1),\ t > 0\,,\\
	& \widehat{h}(-1,t) = H_{-1}\quad \text{and} \quad \widehat{h}(1,t) = H_1\,,& & \qquad\qquad t > 0\,,\\
	& \widehat{h}(x,0) = h_0(x)\,,& & \qquad \qquad x \in (-1,1)\,.
	\end{aligned}
\end{aligned}
\end{equation}

Existence, weak comparison principles, and asymptotic behavior of bounded weak solutions for this type of initial-boundary value problems have been studied in \cite{Diaz2001, DiazThelin1994}.
Note that the analysis (though very thorough) presented
in that two aforementioned papers, is not complete yet and there are many interesting and difficult open problems waiting for their resolution. 
It turns out
from the results and thorough discussions in
\cite{Diaz2001, DiazThelin1994},
that the properties
of the steady states play important role in the theory for the bounded weak solutions to the boundary-initial value problem 
\eqref{eq:math:model:final}. This is one of the motivations for our detailed study
of the steady-states of
\eqref{eq:math:model:final}
presented in this paper.

For the sake of keeping the amount of technicalities to the minimum, we limit ourselves to the case of {\it the same water level in both ditches}, i.e.,
$H\stackrel{\mathrm{def}}{=} H_{-1} = H_1 > 0$. Note that we also assume $H>0$ in this paper, since the case $H=0$ is
essentially different and much more difficult to analyze. Let us observe that
$
\widehat{h}(x,t)\equiv \mathrm{const.}=H
$ is the steady state of \eqref{eq:math:model:final}
for
$f\equiv 0$, and $\widehat{h}(-1,t)=\widehat{h}(1,t)\equiv \mathrm{const.}=H$.
Our aim is to study perturbations of this steady-state solution caused by rain or evaporation constant in time, which is incorporated into the model by considering $f\not\equiv 0$. 
For this reason, we will write the steady-state
solution in the form $\widehat{h}(x,t) = u(x) + H$ and rewrite the equation for steady states of
\eqref{eq:math:model:final}
in terms of unknown function $u$:
\begin{equation}
\label{eq_h_elliptic}
\begin{aligned}
	&
	- 
	\frac{\rm d}{{\rm d} x}
	\left[
	(u(x) + H) \left|\frac{{\rm d} u}{{\rm d} x}(x) \cos\varphi + \sin\varphi \right|^{p - 2}
	\left(\frac{{\rm d} u}{{\rm d} x}(x) \cos\varphi + \sin\varphi\right)
	\right]
	\\
	&
	\begin{aligned}
	&
	= 
	f(x)\,, & \qquad\qquad x \in (-1,1)\,,\\ 
	&
	u(-1) = u(1) = 0\,,&
	\end{aligned}
\end{aligned}
\end{equation}
where in the case of steady-states
we must assume that the function $f$ is independent of the time variable and thus we write $f(x)$ instead of $f(x,t)$. We will use this notation from now on.
Moreover, it is natural to assume that the water level cannot drop below 
the impermeable layer and thus we will consider only solutions that satisfy $u(x)\geq -H$ for $x\in (-1,1)$.

We define weak solution to the problem \eqref{eq_h_elliptic} as follows.

\begin{definition}[Weak solution]
{\rm
Let $p > 1$, $0 < \varphi < \pi/2$, and $f \in L^1(-1,1)$ be given.
By a  weak solution to \eqref{eq_h_elliptic}
we mean a function $u\in W^{1,p}_0(-1,1)$,
$u\geq -H$,
satisfying 
\begin{equation}
\label{weak_formulation}
\begin{aligned}
& \int_{-1}^1
	(u(x) + H) \left|u'(x) \cos\varphi + \sin\varphi \right|^{p - 2}
	\left(u'(x) \cos\varphi + \sin\varphi\right)v'(x)\mathrm{d}x\\
& = 
\int_{-1}^{1}	f(x) v(x)\mathrm{d}x
\end{aligned}
\end{equation}
for all $v\in W^{1,p}_0(-1,1)$, where $u'$ and $v'$ stands for the respective weak derivatives.
}
\end{definition}

\begin{remark}
{\rm
Let us note that any
$u\in W_0^{1,p}(-1,1)$
has a representative 
in $AC[-1,1]$
(space of absolutely continous functions on $[-1,1]$). Thus
every weak solution in our sense is also a~bounded weak solution. So our definition is compatible with results from \cite{Diaz2001, DiazThelin1994}.

Moreover, let us note that the classical derivative of the  
$AC[-1,1]$ representative of $u\in W^{1,p}(-1,1)$ exists
a.e. in $[-1,1]$, belongs to $L^p(-1,1)$, and coincides with the weak derivative of $u$ a.e. in $[-1,1]$ (see, e.g., \cite[Thm.~2.1.4]{Ziemer89}). Therefore, we will use the same notation $u'$ for both the weak and classical derivatives of $u$ with respect to $x$. Also, when we speak about higher regularity of weak solutions such as smoothness, we mean that this regularity refers to the $AC[-1,1]$ representative of $u\in W^{1,p}(-1,1)$.

Finally, let us note that our definition of weak solution is not the most general one. Indeed, it is enough to assume that 
\begin{equation}
\label{wlpprime}
(u(\cdot) + H) \left|u'(\cdot) \cos\varphi + \sin\varphi \right|^{p - 2}
	\left(u'(\cdot) \cos\varphi + \sin\varphi\right)\in L^{p'}(-1,1)
 \end{equation}
 in order that the integral on the left hand side of \eqref{weak_formulation} makes sense.
Indeed, a~function $u\not\in W_0^{1,p}(-1,1)$ satisfying \eqref{wlpprime}
is, e.g.,
$u(x)=|x|^{\alpha}-1$,
if $H=1$ and 
$$0<\left(\frac{p-1}{p}\right)^2<\alpha \leq \frac{p-1}{p}<1\,.$$
However, dealing with such more general weak solutions is quite delicate.  Especially, choice of proper function spaces and finding additional conditions to obtain physically relevant solutions are difficult questions. 
For that reason, we limit ourselves to weak solutions as defined above in this paper.
} 
\end{remark}

\section{Properties of weak solution}
\label{sec:prop:w:sol}

\subsection{Regularity results.}
In this section we will show that any weak solution $u$ to \eqref{weak_formulation} satisfies $u > - H$ and $u \in C^1[-1,1]$, provided the following hypothesis on $f$ is imposed.
\par
\medskip
\noindent
{\bf Hypothesis (HF).} 
Assume that $f \in L^1(-1,1)$ and satisfies the inequality
\begin{equation*}
%\label{cond:not:touching:ground}
\min\limits_{x_0\in [-1,1]}
\min\limits_{x\in [x_0,1]}
\left(
H^{\frac{p}{p-1}}
+
\frac{p}{(p-1)\cos\varphi}
\int_x^1 \Phi_{p'}\left(\int_{x_0}^{\tau} f(s)\mathrm{d}s
\right)\mathrm{d}\tau
\right)
>0\,,
\end{equation*}
where 
$\Phi_{q}(z)\stackrel{\mathrm{def}}{=}
|z|^{q-1}
\mathop{\mathrm{sign}} z
$, for $q>1$.

\begin{remark}
\label{rem:hf}
{\rm 
Let us note that Hypothesis (HF) is satisfied for any
$H>0$ and $f\geq 0$ a.e. in $(-1,1)$. This will turn out to be very usefull in our results concerning maximum principles 
(Proposition~\ref{wmp} and Theorem~\ref{smp}).}
\end{remark}

Our proofs in this section rely on rewriting the  weak formulation~\eqref{weak_formulation} as a differential 
equation of the first order.
\par
\medskip
\noindent
\paragraph{\bf From the weak formulation to a differential 
equation of the first order.}
For any test function $v\in C_{\mathrm{c}}^{\infty}(-1,1)$, we can integrate the right-hand side 
of weak formulation~\eqref{weak_formulation} by parts and obtain
\begin{multline*}
%\label{weak_formulation_integrated}
\int_{-1}^1
	\left(u(x) + H\right) \left|u'(x) \cos\varphi + \sin\varphi \right|^{p - 2}
	\left(u'(x) \cos\varphi + \sin\varphi\right)v'(x)\mathrm{d}x
= \\
-\left[v(x)\int_x^1 f(\sigma)\mathrm{d}\sigma\right]_{-1}^{1} + \int_{-1}^{1}  v'(x)\int_x^1 f(\sigma)\mathrm{d}\sigma \mathrm{d}x
\end{multline*}
which yields
\begin{multline}
\label{weak_formulation_integrated}
\int_{-1}^1
\left[
	\left(u(x) + H\right) \left|u'(x) \cos\varphi + \sin\varphi \right|^{p - 2}
	\left(u'(x) \cos\varphi + \sin\varphi\right) \right.\\
	\left. - \int_x^1 f(\sigma)\mathrm{d}\sigma
\right]	
	v'(x)\mathrm{d}x
= 0
\end{multline}
for all $v\in C_{\mathrm{c}}^{\infty}(-1,1)$.
This means that the distributional derivative of the expression ``$[\dots]$'' in \eqref{weak_formulation_integrated} is zero.
Taking into account that $u\in  W^{1,p}_0(-1,1)\hookrightarrow
AC[-1,1]$, the expression 
``$[\dots]$'' belongs to $L^{1}(-1,1)$. Thus by \cite[Lem. 1.2.1., p.~13]{JostLiJost1998}, we obtain
\begin{multline}
\label{eq_der_const_+}
\left(u(x) + H\right) \left|u'(x) \cos\varphi + \sin\varphi \right|^{p - 2}
	\left(u'(x) \cos\varphi+{} \right.\\
	\left. {}+ \sin\varphi\right) - \int_x^1 f(\sigma)\mathrm{d}\sigma = \\ 
	\kappa\equiv \mathrm{const.}\quad\mbox{a.e. in } [-1,1]\,.
\end{multline}

\begin{theorem}
    \label{thm:regularity:restriction}
     Let $u$ be a weak solution to \eqref{eq_h_elliptic} and $x_{1} \in [-1, 1)$ be such that 
        \begin{equation}
        \label{easy_est_from_below_delta}
            u(x) > - H
            \mbox{ for all } [x_{1}, 1]\,.
        \end{equation}
    Then $u|_{[x_{1}, 1]} \in C^1[x_{1}, 1]$.
\end{theorem}

\paragraph{\em Proof.}
The assumption \eqref{easy_est_from_below_delta}  allows us to rewrite \eqref{eq_der_const_+} as
\begin{multline}
\label{eq_der_gH}
 \left|u'(x) \cos\varphi + \sin\varphi \right|^{p - 2}
	\left(u'(x) \cos\varphi
	 + \sin\varphi\right)  = \\
	 \frac{1}{\left(u(x) + H\right)}
	 \left(
	 \kappa + \int_x^1 f(\sigma)\mathrm{d}\sigma
	 \right)
	 \quad\mbox{a.e. in } [x_{1},1]\,.
\end{multline}
Then we obtain:
\begin{equation}
\label{eq:uprime:ae:kappa}
u'(x)=
\frac{1}{\cos\varphi}
\,
\Phi_{p'}\left(
\frac{1}{\left(u(x) + H\right)}
	 \left(
	 \kappa + \int_x^1 f(\sigma)\mathrm{d}\sigma
	 \right)
\right) - \tan\varphi
	 \quad\mbox{a.e. in } [x_{1},1]\,.
\end{equation}
In the remainder of this proof, $u$ will be identified with the $AC[-1,1]$ representative of $u\in W_0^{1,p}(-1,1)$.  
Now let us define  function
$F \colon \mathbb{R} \times [x_{1},1] \to\mathbb{R}$,
$$
F(z,x)\stackrel{
\mathrm{def}}{=}
z-
\frac{1}{\cos\varphi}
\,
\Phi_{p'}\left(
\frac{1}{\left(u(x) + H\right)}
	 \left(
	 \kappa + \int_x^1 f(\sigma)\mathrm{d}\sigma
	 \right)
\right) + \tan\varphi\,.
$$
It follows 
from 
%$u\in W_0^{1,p}
%(-1,1)\hookrightarrow %AC[-1,1]$
$u\in AC[-1,1]$
and~\eqref{easy_est_from_below_delta} that the function $F$
is continuous on $\mathbb{R}\times [x_{1},1]$.
Moreover, it is strictly increasing in the first variable and
$$
\lim_{z\to\pm\infty} F(z,x)
=\pm\infty\mbox{ for any }
x\in[x_{1},1]\,.
$$
Thus for any $x\in [x_{1},1]$
there is unique $z(x)$ such that
\begin{equation}
\label{fzx}
F(z(x), x) = 0\,.
\end{equation}
Indeed,
$$
    z(x) = \frac{1}{\cos\varphi}
\,
\Phi_{p'}\left(
\frac{1}{\left(u(x) + H\right)}
	 \left(
	 \kappa + \int_x^1 f(\sigma)\mathrm{d}\sigma
	 \right)
\right) + \tan\varphi\,,
$$
which is continuous on $[x_1,1]$ since $u \in AC[-1,1]$ and 
\eqref{easy_est_from_below_delta} holds.
From~\eqref{eq:uprime:ae:kappa}, we also know that
$F(u'(x),x)=0$ a.e. in $[x_{1},1]$. Since $z(x)$ is the only solution to~\eqref{fzx} for any given $x\in [x_{1},1]$, it must hold
$$
z(x)=u'(x) \mbox{ a.e. in }[x_{1},1]\,.
$$
Now, as $u\in AC[-1,1]$, we have
$u(x)=\int_{x_1}^x u'(\sigma)\mathrm{d}\sigma =
\int_{x_1}^x z(\sigma)\,\mathrm{d}\sigma$ for any $x\in [x_1, 1]$. Taking into account that $z\in C[x_1, 1]$, we see that the classical derivative of
$u$ exists in every point
in $x\in (x_1, x)$ and it is continuously extendable to $[x_1,1]$. Hence the restriction to $[x_{1},1]$ of $u\in AC[-1,1]$ 
belongs to $C^{1}[x_{1},1]$. 
\hfill$\blacksquare$
\par\medskip

\begin{corollary}
\label{cor:first:order:formula}
    Let $u$ be a weak solution to \eqref{eq_h_elliptic} and $x_1$ is as in Theorem~\ref{thm:regularity:restriction}. Then
    \begin{multline}
        \label{eq_der_kappa}
        \left(u(x) + H\right) \left|u'(x) \cos\varphi + \sin\varphi \right|^{p - 2}
	   \left(u'(x) \cos\varphi + \sin\varphi\right)
        - 
        \\ 
        H \left|u'(1) \cos\varphi + \sin\varphi \right|^{p - 2} 
	   \left(u'(1) \cos\varphi + \sin\varphi\right) = 
        \\
        \int_x^1 f(\sigma)\mathrm{d}\sigma\,
\end{multline}
holds pointwise everywhere in $[x_1,1]$.
\end{corollary}
\paragraph{\em Proof.}
Since $u|_{[x_{1}, 1]} \in C^1[x_{1}, 1]$ by Theorem~\ref{thm:regularity:restriction}, we may evaluate \eqref{eq_der_const_+} at $x = 1$ to obtain
\begin{multline}
\label{expr_kappa}
\kappa = \\ (u(1) + H) \left|u'(1) \cos\varphi + \sin\varphi \right|^{p - 2}
	\left(u'(1) \cos\varphi + \sin\varphi\right) - \int_{1}^{1} f(\sigma)\mathrm{d}\sigma = \\
	 H \left|u'(1) \cos\varphi + \sin\varphi \right|^{p - 2}
	\left(u'(1) \cos\varphi + \sin\varphi\right)\,.
\end{multline}
Combining \eqref{eq_der_const_+} and \eqref{expr_kappa},
we obtain \eqref{eq_der_kappa}
which holds pointwise everywhere in $[x_1,1]$.
\hfill$\blacksquare$
\par\medskip

Now, assuming Hypothesis {\rm (HF)},
we will show 
that $u > - H$ on the whole $[-1, 1]$ and, 
consequently, $u \in C^1[-1,1]$.

\begin{lemma}
\label{lem:est:alpha}
Let $x_0\in(-1, 1)$ and
$u$ be a weak solution to \eqref{eq_h_elliptic}
such that
$u(x_0)=-H$ 
and 
$u(x)>-H$ for all 
$x\in (x_0,1)$.
Then
\begin{equation}
\label{alpha:est}
 H 
 \left|u'(1) \cos\varphi + \sin\varphi \right|^{p - 2}
	\left(u'(1) \cos\varphi + \sin\varphi\right)
 = 
- \int_{x_0}^1 f(x)\mathrm{d}x\,.
\end{equation}
\end{lemma}

\paragraph{\em Proof.}
Let us use
$$
\displaystyle
v_{\varepsilon}(x)
\stackrel{\mathrm{def}}{=}
\left\{
\begin{array}{rl}
0 &\mbox{for } x\in (-1, x_0] \,, \\[0.15cm]
 (x-x_0)/\varepsilon &\mbox{for } x\in (x_0, x_0+\varepsilon)\,, \\[0.15cm]
 1   &\mbox{for } x\in [x_0+\varepsilon, 1-\varepsilon]\,, \\[0.15cm]
 (1-x)/\varepsilon &\mbox{for } x\in (1-\varepsilon, 1)
\end{array}
\right.
$$
as a test function in \eqref{weak_formulation}, with $0<\varepsilon<(1-x_0)/2$.
Then 
\begin{multline*}
\int_{-1}^1(u(x)+H)
\left|
u'(x)\cos\varphi + \sin\varphi
\right|^{p-2}
\left(
u'(x)\cos\varphi + \sin\varphi
\right)v_{\varepsilon}'(x)\mathrm{d}x 
\\[0.2cm]=
\int_{-1}^1 f(x)v_{\varepsilon}(x)\mathrm{d}x
\end{multline*}
becomes
\begin{multline}
\frac{1}{\varepsilon}
\int_{x_0}^{x_0+\varepsilon}
(u(x)+H)
\left|
u'(x)\cos\varphi + \sin\varphi
\right|^{p-2}
\left(
u'(x)\cos\varphi + \sin\varphi
\right)\mathrm{d}x 
\\[0.2cm] -
\frac{1}{\varepsilon}
\int_{1-\varepsilon}^1
(u(x)+H)
\left|
u'(x)\cos\varphi + \sin\varphi
\right|^{p-2}
\left(
u'(x)\cos\varphi + \sin\varphi
\right)\mathrm{d}x
\\[0.2cm]
=
\int_{-1}^1 f(x)v_{\varepsilon}(x)\mathrm{d}x\,.
\end{multline}
Using the fact that $u\in C^1[x_0+\delta, 1]$ for any $0 < \delta < 1 - x_0$ by Theorem~\ref{thm:regularity:restriction},
we find that
\begin{multline}
\lim\limits_{\varepsilon\to 0+}
\frac{1}{\varepsilon}
\int_{1-\varepsilon}^{1}
(u(x)+H)
\left|
u'(x)\cos\varphi + \sin\varphi
\right|^{p-2}
\left(
u'(x)\cos\varphi + \sin\varphi
\right)\mathrm{d}x \\[0.2cm] =
H \left|
u'(1)\cos\varphi + \sin\varphi
\right|^{p-2}
\left(
u'(1)\cos\varphi + \sin\varphi
\right)
\end{multline}
since $x=1$ is a Lebesgue point of the integrand and $u(1)=0$.
By the Lebesgue dominated convergence theorem, we also find that
$$
\lim\limits_{\varepsilon\to 0+}
\int_{-1}^1 f(x)v_{\varepsilon}(x)\mathrm{d}x
=
\int_{x_0}^1 f(x)\mathrm{d}x\,.
$$
Thus
\begin{multline}
\label{pom:lim}
 H 
 \left|u'(1) \cos\varphi + \sin\varphi \right|^{p - 2}
	\left(u'(1) \cos\varphi + \sin\varphi\right)
 = 
- \int_{x_0}^1 f(x)\mathrm{d}x \\
+
\lim\limits_{\varepsilon\to 0+}
\frac{1}{\varepsilon}
\int_{x_0}^{x_0+\varepsilon}
(u(x)+H)
\left|
u'(x)\cos\varphi + \sin\varphi
\right|^{p-2}
\left(
u'(x)\cos\varphi + \sin\varphi
\right)\mathrm{d}x
\,.
\end{multline}
We will prove \eqref{alpha:est}
by showing that the last term vanishes. 
For any $x\in (x_0, x_0+\varepsilon)$,
we have
\begin{multline}
\label{pom:1}
\frac{|u(x)+H|}{\varepsilon^{1/p'}}
=
\frac{|u(x)-(-H)|}{\varepsilon^{1/p'}}=
\frac{
\left|
\int_{x_0}^x u'(s)\mathrm{d}s
\right|
}{\varepsilon^{1/p'}}
\leq
\frac{
\int_{x_0}^x
\left|u'(s)\right|\mathrm{d}s
}{(x-x_0)^{1/p'}}\\[0.2cm]
\leq
\frac{
\left(
\int_{x_0}^x 1^{p'}\mathrm{d}s
\right)^{1/p'}
}{(x-x_0)^{1/p'}}
\left(
\int_{x_0}^x |u'(s)|^p \mathrm{d}s
\right)^{1/p}
\leq \left(
\int_{x_0}^{x_0+\varepsilon} |u'(s)|^p \mathrm{d}s
\right)^{1/p}\,.
\end{multline}

On the other hand,
\begin{multline}
\label{pom:2}
\frac{
\int_{x_0}^{x_0+\varepsilon}
|u'(x)\cos\varphi + \sin\varphi|^{p-1}\mathrm{d}x}{\varepsilon^{1/p}} \\
\leq
\frac{
\left(
\int_{x_0}^{x_0+\varepsilon}
1^p\mathrm{d}x
\right)^{1/p}
}
{\varepsilon^{1/p}}
\left(\int_{x_0}^{x_0+\varepsilon}
|u'(x)\cos\varphi + \sin\varphi|^p\mathrm{d}x\right)^{1/p'}
\\ =
\left(\int_{x_0}^{x_0+\varepsilon}
|u'(x)\cos\varphi + \sin\varphi|^p\mathrm{d}x\right)^{1/p'}\,.
\end{multline}
Now, combining \eqref{pom:1} and \eqref{pom:2}, we obtain
\begin{multline}
\lim\limits_{\varepsilon\to 0+}
\frac{1}{\varepsilon}
\left|\int_{x_0}^{x_0+\varepsilon}
(u(x)+H)
\left|
u'(x)\cos\varphi + \sin\varphi
\right|^{p-2}
\left(
u'(x)\cos\varphi + \sin\varphi
\right)\mathrm{d}x\right| 
\\ \leq
\lim\limits_{\varepsilon\to 0+}
\frac{1}{\varepsilon^{1/p}}
\int_{x_0}^{x_0+\varepsilon}
\frac{|u(x)+H|}{\varepsilon^{1/p'}}
\left|
u'(x)\cos\varphi + \sin\varphi
\right|^{p-1}
\mathrm{d}x \\ \leq
\lim\limits_{\varepsilon\to 0+}
\frac{1}{\varepsilon^{1/p}}
\int_{x_0}^{x_0+\varepsilon}
\left(
\int_{x_0}^{x_0+\varepsilon} |u'(x)|^p \mathrm{d}x
\right)^{1/p}
\left|
u'(x)\cos\varphi + \sin\varphi
\right|^{p-1}
\mathrm{d}x \\ \leq
\lim\limits_{\varepsilon\to 0+}
\left(
\int_{x_0}^{x_0+\varepsilon} |u'(x)|^p \mathrm{d}x
\right)^{1/p}
\left(\int_{x_0}^{x_0+\varepsilon}
|u'(x)\cos\varphi + \sin\varphi|^p\mathrm{d}x\right)^{1/p'}= 0\,,
\end{multline}
since $u\in W_0^{1,p}(-1,1)$.
Thus the last term in \eqref{pom:lim}
vanishes and the proof is complete.
\hfill$\blacksquare$
\par\medskip

\begin{theorem}
\label{thm:no:touch}
Let Hypothesis {\rm (HF)} hold. Then any weak solution $u$ to \eqref{eq_h_elliptic}
satisfies $u>-H$ on $[-1,1]$.
\end{theorem}

\paragraph{\em Proof.}
Assume by contradiction that there is at least one $\xi\in (-1,1)$ such that
$u(\xi)=-H$. Since $u\in W_0^{1,p}(-1,1)\hookrightarrow AC[-1,1]$ and $u(1)=0$, we may choose $\xi = x_0$ 
such that $u(x)>-H$ for all $x\in (x_0, 1)$.
Then by Corollary~\ref{cor:first:order:formula}, we have
\begin{multline}
%\label{eq_der_kappa}
 H 
 \left|u'(1) \cos\varphi + \sin\varphi \right|^{p - 2}
	\left(u'(1) \cos\varphi + \sin\varphi\right)
- 
\\ \left(u(x) + H\right) \left|u'(x) \cos\varphi + \sin\varphi \right|^{p - 2} 
	\left(u'(x) \cos\varphi + \sin\varphi\right) = \\
- \int_{x}^1 f(\sigma)\mathrm{d}\sigma
\end{multline}
for all $x \in [x_0 + \delta, 1)$ with any $0 < \delta < 1 - x_0$.
Using \eqref{alpha:est}, we find 
\begin{multline}
-\left(u(x) + H\right) \left|u'(x) \cos\varphi + \sin\varphi \right|^{p - 2} 
	\left(u'(x) \cos\varphi + \sin\varphi\right) = \\
- \int_{x}^1 f(\sigma)\mathrm{d}\sigma
+ \int_{x_0}^1 f(\sigma)\mathrm{d}\sigma \\
=
\int_{x_0}^x f(\sigma)\mathrm{d}\sigma
\end{multline}
for any $x\in (x_0,1)$.
Using substitution
$w(x) = u(x)+H$ and denoting $\tilde F(x)=\int_{x_0}^x f(s)\mathrm{d}s$, we obtain
$$
w(x)
\Phi_{p}\left(
w'(x)
\cos\varphi+
\sin \varphi
\right)
=
- \tilde F(x)
$$
and, equivalently,
$$
\Phi_{p}\left([w(x)]^{\frac{1}{p-1}}w'(x)
\cos\varphi+
[w(x)]^{\frac{1}{p-1}} \sin \varphi
\right)
=
- \tilde F(x)\,.
$$
Since $\Phi_{p'}$ is an inverse to $\Phi_p$, we get
$$
[w(x)]^{\frac{1}{p-1}}w'(x)
\cos\varphi+
[w(x)]^{\frac{1}{p-1}} \sin \varphi=
-\Phi_{p'}\left(\tilde F(x)
\right)
$$
Using substitution
$v(x) = [w(x)]^{\frac{p}{p-1}}$,
we have the following ODE
$$
\begin{aligned}
\frac{p-1}{p} v^{\prime}(x) \cos\varphi+
v^{1 / p}(x) \sin \varphi & =
- \Phi_{p'}\left(
\tilde F(x)
\right)\,,\\
v(1)& = H^{\frac{p}{p-1}}\,,\\
v(x) & \geq 0\,,\\
\varphi & > 0\,.
\end{aligned}
$$
Since $v(x) \geq 0$, we have
$$
\frac{p-1}{p} v^{\prime}(x) \cos\varphi
\leq
- \Phi_{p'}\left(
\tilde F(x)
\right)
$$
and integrating from $x$ to 1 we get
$$
\frac{(p-1)\cos\varphi}{p} \left(v(1)-v(x)\right) 
\leq
-\int_x^1 
 \Phi_{p'}\left(\tilde F(\tau)\right)\,\mathrm{d}\tau\,.
$$
It follows that
$$
v(x)
\geq
H^{\frac{p}{p-1}}+\frac{p}{(p-1)\cos\varphi}\int_x^1 \Phi_{p'}\left(\tilde F(\tau)\right)\,\mathrm{d}\tau\,.
$$
By  (HF),
we have
$$
K=\min\limits_{x_0\in [-1,1]}\min\limits_{x\in [x_0,1]}
\left(
H^{\frac{p}{p-1}}
+
\frac{p}{(p-1)\cos\varphi}
\int_x^1 \Phi_{p'}\left(\int_{x_0}^{\tau} f(s)\mathrm{d}s\right)\,\mathrm{d}\tau
\right) > 0\,,
$$
then $v(x) \geq K$ must hold for all $x\in (x_0,1)$. This contradicts $v(x_0)=0$
as $v$ is continuous on $[-1,1]$.
\hfill$\blacksquare$
\par\medskip

\begin{theorem}[$C^1$-regularity of a weak solution]
\label{thm:C1:reg}
    Let Hypothesis {\rm (HF)} hold. Then any weak solution to \eqref{eq_h_elliptic} satisfies $u \in C^1[-1,1]$ .
\end{theorem}
\paragraph{\em Proof.} By Theorem~\ref{thm:no:touch}, $u(x) > -H$ for all $x \in [-1,1]$. 
Then we may choose $x_1 = -1$ in Theorem~\ref{thm:regularity:restriction} to obtain that $u \in C^1[-1,1]$.
\hfill$\blacksquare$
\par\medskip

\begin{corollary}
\label{cor:first:order:plus}
Under Hypothesis {\rm (HF)},
the following equation
 \begin{multline}
        \label{eq:first:order:plus}
        \left(u(x) + H\right) \left|u'(x) \cos\varphi + \sin\varphi \right|^{p - 2}
	   \left(u'(x) \cos\varphi + \sin\varphi\right)
        - 
        \\ 
        H \left|u'(1) \cos\varphi + \sin\varphi \right|^{p - 2} 
	   \left(u'(1) \cos\varphi + \sin\varphi\right) = 
        \\
        \int_x^1 f(\sigma)\mathrm{d}\sigma\
\end{multline}
holds pointwise for all $x \in [-1,1]$.
\end{corollary}

\paragraph{\em Proof.} 
The proof of \eqref{eq:first:order:plus} is similar as the proof of Corollary~\ref{cor:first:order:formula} with $x_1 = -1$
since $u \in C^1[-1,1]$ by Theorem~\ref{thm:C1:reg}.
\hfill$\blacksquare$
\par\medskip

\subsection{A~priori bound on a weak solution.}
Here, we use the regularity result from the previous section to obtain a~priori bound on the $L^{\infty}$-norm of a weak solution to \eqref{weak_formulation} depending on $\varphi$ and $\|f\|_{L^1(-1,1)}$. The proof relies on the first order formula \eqref{eq:first:order:plus}.

\begin{theorem}
\label{thm:boundness:h}
Assume that Hypothesis {\rm (HF)} holds. Then any weak solution $u$ to \eqref{eq_h_elliptic}
satisfies
\begin{equation}
\label{est:l1}
\|u\|_{\infty}
\leq 
\frac{\|f\|_{L^1(-1,1)}}{\left(\sin\varphi \right)^{p - 1}}\,.
\end{equation}
\end{theorem}

\paragraph{\em Proof.} We distinguish four basic cases further possibly divided into subcases.

\par\noindent{\it Case 1. $u\equiv 0$.} Statement is satisfied trivially.
 
\par\noindent{\it Case 2: $u(x)\geq 0$, $u(x_{\mathrm{max}})>0$.} Maximum is achieved in the interior point $x_{\mathrm{max}}\in (-1,1)$. Thus $u'(x_{\mathrm{max}})=0$ since $u \in C^1[-1,1]$. Minimum is achieved at the boundary.
Evaluating \eqref{eq:first:order:plus} at $x_{\mathrm{max}}$, we get
\begin{multline}
\label{int_from_max}
\left(u(x_{\mathrm{max}})+ H\right) \left(\sin\varphi \right)^{p - 1}
- 
\\ H \left|u'(1) \cos\varphi + \sin\varphi \right|^{p - 2} 
	\left(u'(1) \cos\varphi + \sin\varphi\right) = \\
\int_{x_{\mathrm{max}}}^1 f(\sigma)\mathrm{d}\sigma\,.
\end{multline}
Since $u\in C^{1}[-1,1]$, 
$u(-1)=0=u(1)$, and $u(x)\geq 0$, we have 
$$
    u'(-1)\geq 0 \geq u'(1)\,.
$$

\par\noindent{\it Subcase 2a: $u'(1) < -\tan\varphi$.} Then $u'(1) \cos\varphi + \sin\varphi<0$. 
Hence 
$$
-  H \left|u'(1) \cos\varphi + \sin\varphi \right|^{p - 2} 
	\left(u'(1) \cos\varphi + \sin\varphi\right) 
	=
	 H \left|u'(1) \cos\varphi + \sin\varphi \right|^{p - 1}
	 > 0\,.
$$
Then it follows from \eqref{int_from_max} that
\begin{multline*}
   0 < u(x_{\mathrm{max}}) \leq u(x_{\mathrm{max}})+ H
+
\frac{H}{\left(\sin\varphi \right)^{p - 1}}
\left|u'(1) \cos\varphi + \sin\varphi \right|^{p - 1}
=
\\
 \frac{1}{\left(\sin\varphi \right)^{p - 1}}\int_{x_{\mathrm{max}}}^1 f(\sigma)\mathrm{d}\sigma
\leq \frac{\|f\|_{L^1(-1,1)}}{\left(\sin\varphi \right)^{p - 1}}\,.
\end{multline*}

\par\noindent{\it Subcase  2b: $0 \geq u'(1) \geq - \tan\varphi$.} Then $0 \leq u'(1) \cos\varphi + \sin\varphi \leq \sin\varphi$. Thus it follows from \eqref{int_from_max} that
\begin{multline*}
0 < u(x_{\mathrm{max}}) \leq \\
 u(x_{\mathrm{max}}) + H\left(1 -  \left(\frac{\left|u'(1) \cos\varphi + \sin\varphi \right|}{\sin\varphi}\right)^{p - 1}  
\right) 
=
\\
\frac{1}{\left(\sin\varphi \right)^{p - 1}}\int_{x_{\mathrm{max}}}^1 f(\sigma)\mathrm{d}\sigma
\leq \frac{\|f\|_{L^1(-1,1)}}{\left(\sin\varphi \right)^{p - 1}}\,.
\end{multline*}

\par\noindent{\it Case 3: $u(x)\leq 0$, $u(x_{\mathrm{min}})<0$.} Minimum is achieved in the interior point $x_{\mathrm{min}}\in (-1,1)$. Thus $u'(x_{\mathrm{min}})=0$ since $u \in C^1[-1,1]$. Maximum is achieved at the boundary. 
From $u\in C^{1}[-1,1]$, 
$u(-1)=0=u(1)$, and $u(x)\leq 0$, we conclude 
$$
    u'(-1)\leq 0 \leq u'(1)\,.
$$
Evaluating \eqref{eq:first:order:plus} at $x_{\mathrm{min}}$, we get
\begin{multline}
\label{int_from_min}
\left(u(x_{\mathrm{min}})+ H\right) \left(\sin\varphi \right)^{p - 1}
- 
\\ H \left|u'(1) \cos\varphi + \sin\varphi \right|^{p - 2} 
	\left(u'(1) \cos\varphi + \sin\varphi\right) = \\
 \int_{x_{\mathrm{min}}}^1 f(\sigma)\mathrm{d}\sigma\,.
\end{multline}
Since $u(x_{\mathrm{min}}) < 0$ and $u'(1) > 0$, we have 
$$
\left(|u(x_{\mathrm{min}})| - H\right) \left(\sin\varphi \right)^{p - 1}
+ 
H \left|u'(1) \cos\varphi + \sin\varphi \right|^{p - 1} 
=
- \int_{x_{\mathrm{min}}}^1 f(\sigma)\mathrm{d}\sigma\,,
$$
which follows
\begin{multline*}   
0 < |u(x_{\mathrm{min}})| \leq |u(x_{\mathrm{min}})| 
+ 
H \left(\left(\frac{\left|u'(1) \cos\varphi + \sin\varphi \right|}{\sin\varphi}\right)^{p - 1} - 1
\right)
= 
\\
-\frac{1}{\left(\sin\varphi \right)^{p - 1}}\int_{x_{\mathrm{max}}}^1 f(\sigma)\mathrm{d}\sigma
\leq \frac{\|f\|_{L^1(-1,1)}}{\left(\sin\varphi \right)^{p - 1}}\,.
\end{multline*}

\par\noindent{\it Case 4: $u(x_{\mathrm{min}})<0<u(x_{\mathrm{max}})$.} 
In this case, 
$x_{\mathrm{min}}\not=x_{\mathrm{max}}$ and $x_{\mathrm{min}},x_{\mathrm{max}}\in (-1,1)$. Without loss of generality assume that $x_{\mathrm{min}} < x_{\mathrm{max}}$.
Since $u \in C^1[-1,1]$, $u'(x_{\mathrm{min}})=0=u'(x_{\mathrm{max}})$. Hence both \eqref{int_from_max} and \eqref{int_from_min} are valid.
Subtracting \eqref{int_from_min} from \eqref{int_from_max}, we obtain
$$
0 < \left(u(x_{\mathrm{max}})-u(x_{\mathrm{min}})\right)\left(\sin\varphi \right)^{p - 1}\leq -\int_{x_{\mathrm{min}}}^{x_{\mathrm{max}}} f(\sigma)\mathrm{d}\sigma
\leq \|f\|_{L^1(-1,1)}\,.
$$

Since $u(x_{\mathrm{min}})<0$, we get
$$
0 < u(x_{\mathrm{max}})\leq u(x_{\mathrm{max}}) - u(x_{\mathrm{min}}) \leq \frac{\|f\|_{L^1(-1,1)}}{\left(\sin\varphi \right)^{p - 1}}\,.
$$
Moreover, by 
$u(x_{\mathrm{max}})>0$, we also get
$$
0 <|u(x_{\mathrm{min}})| < u(x_{\mathrm{max}}) - u(x_{\mathrm{min}}) \leq \frac{\|f\|_{L^1(-1,1)}}{\left(\sin\varphi \right)^{p - 1}}\,.$$

Finally, we find that 
$$
   0 \leq \|u\|_{\infty}\leq \frac{\|f\|_{L^1(-1,1)}}{\left(\sin\varphi \right)^{p - 1}}
$$
since these four are the all possible cases.
\hfill$\blacksquare$
\par\medskip

\subsection{Existence of weak solution}
Now we establish an existence result for the problem \eqref{eq_h_elliptic} using combination of a~priori estimate \eqref{est:l1} and classical theory of pseudomonotone operators.
Since this theory is well known, we refer the reader 
to, e.g.,  \cite[p.~5; Def.~2.1, pp.~31--32; and Def.~2.5, p.~33]{Roubicek2ndedit} for definitions of bounded, coercive, and pseudomonotone operator, respectively. As usual,
dual of $W_0^{1,p}(-1,1)$ will be denoted by $W^{-1,p'}(-1,1)$.
Let us note that a~priori estimate \eqref{est:l1} plays a key role in verification
that the nonlinear term in \eqref{eq_h_elliptic}
satisfies structural conditions (presented in \cite{Roubicek2ndedit}),
which allow application of the theory of pseudomonotone operators. Our existence result reads as follows.  
\begin{comment}
Now, we use the estimate on $u$ to obtain an existence result for the problem \eqref{eq_h_elliptic}. 
The proof of the result, see Theorem~\ref{thm:existence} below,
follows the standard steps as, e.g., in Roub\'i\v{c}ek~\cite[Section~2.4, pp.~42--55]{Roubicek2010} and hence the technical part of the proof is postponed to the Appendix A. 
We refer the reader to \cite[p.~2]{Roubicek2010}, \cite[Definition~2.1, pp.~31--32]{Roubicek2010}, and \cite[Definition~2.5, p.~33]{Roubicek2010} for the definition of bounded, coercive, and monotone operator, respectively.
\end{comment}
%{\color{red} V nasledujici vete pisu $p>1$, dukaz lze snadno upravit, pridam to pozdeji do appendixu.
%https://math.stackexchange.com/questions/499503/how-to-show-that-p-laplacian-operator-is-monotone
%}
\begin{theorem}[Existence of weak solution]
\label{thm:existence}
Assume that
\begin{equation}
\label{small:f}
   \|f\|_{L^1(-1,1)} < H(\sin\varphi)^{p - 1}\,.
\end{equation}
Then the problem \eqref{eq_h_elliptic} possesses at least one weak solution.
\end{theorem}
\paragraph{\em Proof.} 
Set 
$$
    k\stackrel{\rm def}{=} \frac{\|f\|_{L^1(-1,1)}}{\left(\sin\varphi \right)^{p - 1}}
$$
and define truncation function
$$
    T_k(r) \stackrel{\rm def}{=} 
    \left\{
    \begin{aligned}
       -k\,,& 
       \qquad r < - k \\
        r, & \quad
        -k \leq r \leq k\,, \\
        k, & \qquad r > k\,.
    \end{aligned}
    \right.
$$
We say that $w\in W_0^{1,p}(-1,1)$ satisfies truncated version of \eqref{weak_formulation} if
\begin{multline}
\label{eq:trunk}
\int_{-1}^1
	(T_k(w(x))+H) \left|w'(x) \cos\varphi + \sin\varphi \right|^{p - 2}
	\left(w'(x) \cos\varphi + \sin\varphi\right)v'(x)\mathrm{d}x= \\ \int_{-1}^{1}	f(x) v(x)\mathrm{d}x
\end{multline}
for all $v\in W^{1,p}_0(-1,1)$. 
Now we briefly comment on regularity properties and 
a~priori bounds valid for solutions to \eqref{eq:trunk}, omitting details since the proofs are analogous to and even simpler than for the case of solutions to \eqref{weak_formulation}.
Considering \eqref{small:f}, we find that $T_k(w(x))>-H$ for all $x\in [-1,1]$. This ensures that the argument used in the proof of Theorem~\ref{thm:regularity:restriction} can be applied directly to the entire interval $[-1,1]$. Consequently, we obtain $w\in C^1[-1,1]$. An analogous result to Corollary~\ref{cor:first:order:plus} holds for solutions to \eqref{eq:trunk}, without assuming Hypothesis {\rm (HF)}. Specifically, any solution $w$ to \eqref{eq:trunk} satisfies
 \begin{multline}        \label{eq:trunc:first:order:plus}
        \left(T_k(w(x)) + H\right) \left|w'(x) \cos\varphi + \sin\varphi \right|^{p - 2}
	   \left(w'(x) \cos\varphi + \sin\varphi\right)
        - 
        \\ 
        H \left|w'(1) \cos\varphi + \sin\varphi \right|^{p - 2} 	   \left(w'(1) \cos\varphi + \sin\varphi\right) = 
        \\
        \int_x^1 f(\sigma)\mathrm{d}\sigma\
\end{multline}
holds pointwise for all $x \in [-1,1]$.
Finaly, a~priori estimate 
\begin{equation}
\label{eq:w:est}
\|w\|_{\infty}
\leq 
\frac{\|f\|_{L^1(-1,1)}}{\left(\sin\varphi \right)^{p - 1}}
\end{equation}
follows from \eqref{eq:trunc:first:order:plus} using the same steps as in the proof of Theorem~\ref{thm:boundness:h}. Now, in view of \eqref{eq:w:est}, \eqref{small:f}, and the definition of $k$, we see that any solution to \eqref{eq:trunk} also satisfies \eqref{weak_formulation}.

We now proceed to establish the existence of solution to \eqref{eq:trunk}, which will then imply the existence of solution to \eqref{weak_formulation}, which is a weak solutions to \eqref{eq_h_elliptic}.
This will be achieved by
verifying, as detailed in Appendix~\ref{app:a}, that the operator $A\colon W^{1,p}_0(-1,1) \to W^{-1,p'}(-1,1)$ defined by the left-hand side of~\eqref{eq:trunk} is a~bounded, coercive, and pseudomonotone operator.
By~\cite[Thm.~2.6, p.~33]{Roubicek2ndedit},
there exists a solution $w\in W^{1,p}_0(-1,1)$ to the operator equation $A(w) = g$
for all $g \in W^{-1,p'}(-1,1)$ and hence also for $g\in W^{-1,p'}(-1,1)$ defined by the right-hand side of~\eqref{eq:trunk}, 
which concludes the proof.
\hfill$\blacksquare$
\par\medskip

\begin{remark}
{\rm
In Theorem~\ref{thm:existence}, we obtained solutions to \eqref{weak_formulation} using the fact that any solution to \eqref{eq:trunk} is also a solution to \eqref{weak_formulation} under the condition \eqref{small:f}. This raises the question of whether \eqref{weak_formulation} might possess solutions that do not also satisfy \eqref{eq:trunk}. While we can only provide partial answer to this question, we can establish the following: If we assume Hypothesis {\rm (HF)} in addition to~\eqref{small:f}, then any weak solution to \eqref{weak_formulation} must also satisfy \eqref{eq:trunk}. Consequently, the two problems become equivalent under this combined assumption.
}   
\end{remark}

\subsection{Linearization at the trivial solution.}
\label{par:linearizaton}
The method of linearization of the $p$-Laplacian at some given solution has been succesfuly used in treating various questions such as validity of comparison principles for elliptic and parabolic problems involving $p$-Laplacian, see, e.g., \cite{BenediktGirgKotrlaTakac2019, CuestaTakac1998, CuestaTakac2000, GueddaVeron1989}, 
Fredholm Alternative for the $p$-Laplacian \cite{DrabekGirgTakac2004, Takac2002}.
The reader who is further interested in this method is refered to \cite{Takac2002}
for the most detailed description of this method.

The main tool of this method is the zero order Taylor formula for the power function with the remainder in the integral form:
\begin{equation}
\label{eq:odd:powers:equality}
    |a|^{p - 2}a - |b|^{p-2}b
    =
   (p - 1)\int_0^1|b + \theta (a - b)|^{p - 2}\,\mathrm{d}\theta\,
   (a - b)
\end{equation}
To estimate the remainder term from below (the more involved estimate), we make use of the following lemma, which is an alternative to \cite[Lem.~A.1, p.~233]{Takac2002} for one-dimensional case.
Main advantage of our approach is that it is better suited for a particular form of reminder term in calculations below and that it provides simple explicit lower bound.
\begin{lemma}
\label{lemma:lower:bound:int}
Let $p > 2$ and $a \in \mathbb{R}$, then
    $$
    (p - 1)\int_0^1|1 + \theta a|^{p - 2}\,\mathrm{d}\theta
    \geq 1/2
    %\left(1 - %\left(\frac{1}{2}\ri%ght)^{p - 1}\right)\,.
$$
\end{lemma}

\paragraph{\em Proof.}
We distinguish two cases.
\par\noindent{\it Case 1.} Let $a\geq -2$. Then
$1+\theta a \geq 1-2\theta\geq 0$ for 
$\theta\in [0, 1/2]$.
Thus 
\begin{multline*}
(p-1)\int_0^1 
\left|
1+\theta a
\right|^{p-2}
\mathrm{d}
\theta
\geq \\
(p-1)\int_0^{1/2} 
\left|
1+\theta a
\right|^{p-2}
\mathrm{d}
\theta
\geq
(p-1)\int_0^{1/2} 
\left|
1-2\theta
\right|^{p-2}
\mathrm{d}
\theta
= 1/2\,.
\end{multline*}
\par\noindent{\it Case 2.} Let $a<-2$. Then
$1+\theta a \leq 1-2\theta\leq 0$ for 
$\theta\in [1/2, 1]$ and we again obtain
\begin{multline*}
(p-1)\int_0^1 
\left|
1+\theta a
\right|^{p-2}
\mathrm{d}
\theta
\geq \\
(p-1)\int_{1/2}^{1} 
\left|
1+\theta a
\right|^{p-2}
\mathrm{d}
\theta
\geq
(p-1)\int_{1/2}^{1} 
\left|
1-2\theta
\right|^{p-2}
\mathrm{d}
\theta
= 1/2\,.
\end{multline*}
This completes the proof.
\hfill$\blacksquare$
\par\medskip

Now we are ready to state our main result concerning linearization of \eqref{eq_h_elliptic} at the zero solution.
\begin{lemma}
\label{prop:nonlinear:to:linear}
Let $p>2$ and Hypothesis {\rm (HF)}  be satisfied.
Let $u$ be a weak solution to nonlinear problem \eqref{eq_h_elliptic} 
and for such $u$ let us define
\begin{equation}
\label{eq:D}
D(x)\stackrel{\rm def}{=}
    (u(x)+H)(p-1) \int_{0}^{1}
    \left|\sin \varphi+\theta u^{\prime}(x) \cos \varphi\right|^{p-2} 
    \,\mathrm{~d} \theta \cos \varphi
\end{equation}
for every $x\in [-1,1]$.
Then $D(\cdot)\in C[-1,1]$
and, for all $x\in [-1,1]$,
the following estimate holds
\begin{equation}
\label{eq:lower:bound:D}
D(x)\geq
K\left(\varphi, H, p, \|u^-\|_{\infty}\right)
 >0\,.
\end{equation}
Moreover,
$u$ is also the weak solution to the following linear problem
\begin{equation}
\label{eq:linearization:line:prob}
\left\{
\begin{aligned}
 - \frac{\mathrm{d}}{\mathrm{d}x}\left(D(x)\frac{\mathrm{d}u}{\mathrm{d}x}(x) \right)
 - (\sin \varphi)^{p - 1} \frac{\mathrm{d}u}{\mathrm{d}x}(x) & = f(x)\\
  u(-1) & = 0 = u(1)\,,
\end{aligned}
\right.
\end{equation} 
that is, $u$ satisfies
\begin{equation}
\label{eq:weak:linear}
\int_{-1}^{1} 
D(x)\,u'(x) v'(x)\, \mathrm{d} x
+(\sin \varphi)^{p-1} \, \int_{-1}^{1} u(x) v'(x)\, \mathrm{d} x
=
\int_{-1}^{1} f(x) v(x) \mathrm{d} x
\end{equation}
for all $v\in W^{1,p}_0(-1,1)$.
\end{lemma}
\paragraph{\em Proof.}
Let us consider any (but fixed) weak solution 
$u$ 
to 
\eqref{eq_h_elliptic}.
Since $u \in C^1[-1,1]$ by Theorem~\ref{thm:C1:reg} and
$
    u(x) > - H
$
by Theorem~\ref{thm:no:touch}, there exists $M > 0$ depending on $\|u^-\|_{\infty}$ such that
$$
    \min\limits_{x\in [-1,1]} \left(u(x)+H\right)\geq M > 0\,.
$$
Then the function $D(\cdot)$  given by \eqref{eq:D}
is well defined on $[-1,1]$ and 
$D(\cdot)\in C[-1,1]$ for $p > 2$.

Moreover, we have 
\begin{multline}
\label{eq:Dx}
D(x) = 
(u(x)+H)(p-1) \int_{0}^{1}\left|\sin \varphi+\theta u^{\prime}(x) \cos \varphi\right|^{p-2} \mathrm{~d} \theta \cos \varphi\geq 
\\
M
(\sin \varphi)^{p-2}
 (p-1) \int_{0}^{1}\left|1 + \theta u^{\prime}(x) \cot \varphi\right|^{p-2} \mathrm{~d} \theta
 \cos \varphi \geq 
\\
\frac{M}{2}
(\sin \varphi)^{p-2} \cos \varphi
>0\,,
\end{multline}
where we used Lemma~\ref{lemma:lower:bound:int}
to estimate 
$(p-1)\int_0^1\dots\mathrm{d}\theta\geq 1/2$. Thus setting 
$K(\varphi, H, p, \|u^-\|_{\infty})=\frac{M}{2}
(\sin \varphi)^{p-2} \cos \varphi$, we established validity of \eqref{eq:lower:bound:D}.

Again, let us consider any (but fixed) weak solution 
$u$ 
to 
\eqref{eq_h_elliptic}.
Let the function $D(\cdot)$
be constructed from this fixed $u$ by~\eqref{eq:D}.
It remains to show that $u$
satisfy~\eqref{eq:weak:linear}.
\begin{comment}
We obtained linear second order ODE (3.27) which,  through D, depends on u being a fixed weak solution to (2.8). It remains to show that the fixed u also satisfies (3.28), i.e., u is weak solution to (3.27). 

It remains to show that $u$
satisfy \eqref{eq:weak:linear} with $D(x)$ given by \eqref{eq:D} being the function corresponding to the fixed solution $u$.
\end{comment}
To do this, let us observe that $u \equiv 0$ satisfies  \eqref{weak_formulation}
for $f\equiv 0$. Indeed, we have
\begin{equation}
\label{eq:formula:2}
\int_{-1}^{1} H|\sin \varphi|^{p-2} \sin \varphi \cdot v'(x)\, \mathrm{d} x=\int_{-1}^{1} 0 \cdot v(x)\, \mathrm{d} x
\end{equation}
for any 
$v\in W^{1,p}_0(-1,1)$.
Now subtracting \eqref{eq:formula:2}
from \eqref{weak_formulation}, i.e., from
\begin{equation*}
\label{eq:formula:1}
\begin{aligned}
&
\int_{-1}^{1}(u(x)+H)\left|u'(x) \cos \varphi+\sin \varphi\right|^{p-2}\left(u'(x) \cos \varphi+\sin \varphi\right) v'(x)\, \mathrm{d} x\\
&
=
\int_{-1}^{1} f(x) v(x)\, \mathrm{d} x\,,
\end{aligned}
\end{equation*}
we obtain
\begin{equation}
\label{eq:linearization:subtracted:formulas}
\begin{aligned}
&
\int_{-1}^{1}\biggl[(u(x)+H)\left|u'(x) \cos \varphi+\sin \varphi\right|^{p-2}\left(u'(x) \cos \varphi+\sin \varphi\right)%\right.
\\ 
&
%\left.
-H(\sin \varphi)^{p-1}\biggr] v^{\prime}(x)\,\mathrm{d} x
=\int_{-1}^{1} f(x) v(x)\, \mathrm{d} x\,.
\end{aligned}
\end{equation}
By ``adding zero'' 
$0=u(x)(\sin \varphi)^{p-1}-u(x)(\sin \varphi)^{p-1}$
to the term ``$[\dots]$'' above, we obtain
$$
\begin{aligned}
& \bigl[\dots\bigr]=
(\underbrace{u(x)+H)}_{\geqslant \mathrm{const.}>0} \underbrace{\left[\left|u'(x) \cos \varphi+\sin \varphi\right|^{p-2}\left(u'(x) \cos \varphi+\sin \varphi\right)-|\sin \varphi|^{p-2} \sin \varphi\right]}_{(p-1) \int_{0}^{1}\left|\sin \varphi+\theta u'(x) \cos \varphi\right|^{p-2} \mathrm{d} \theta \, u'(x) \cos \varphi}
\\
&
+u(x)(\sin \varphi)^{p-1} \\
&
= (u(x)+H)(p-1)
\int_{0}^{1}\left|\sin \varphi+\theta u'(x) \cos \varphi\right|^{p-2} \mathrm{d} \theta\,u'(x)\cos\varphi + u(x)(\sin \varphi)^{p-1} = \\
&
D(x)u'(x)+u(x)(\sin \varphi)^{p-1}\,.
\end{aligned}
$$
Using the last formula on the last line above instead of $[\dots]$ in   
\eqref{eq:linearization:subtracted:formulas},
we obtain
$$
\int_{-1}^{1} D(x)\,u'(x) v'(x)\, \mathrm{d} x
+(\sin \varphi)^{p-1} \, \int_{-1}^{1} u(x) v'(x)\, \mathrm{d} x
=
\int_{-1}^{1} f(x) v(x) \mathrm{d} x
$$
for any $v\in W_0^{1,p}(-1,1)$.
This means that $u$ is the weak solution to the linear problem 
\eqref{eq:linearization:line:prob}. This ends the proof.\hfill$\blacksquare$
\par\medskip

Let us note that,
due to \eqref{eq:D}, the diffusion coefficient $D(\cdot)$ depends on a weak solution $u$ under consideration of the nonlinear problem \eqref{eq_h_elliptic}. The following lemma provides lower bound on the diffusion coefficient independent
of a particular choice of a weak solution $u$ to the nonlinear problem \eqref{eq_h_elliptic}.

\begin{lemma}
\label{prop:uniform:lower:bound}
Let $p>2$ and Hypothesis {\rm (HF)}  be satisfied.
Moreover, let
$$
   \|f\|_{L^1(-1,1)} < H(\sin\varphi)^{p - 1}\,.
$$
Then 
\begin{equation}
\label{estimate:Dx}
\begin{aligned}
& D(x) \geq K'\left( \varphi, H, p, \|f\|_{L^1(-1,1)}\right)
\stackrel{\mathrm{def}}{=} \\
& \frac{1}{2}\left(
        H - \frac{1}{(\sin\varphi)^{p - 1}} \|f\|_{L^1(-1,1)}
    \right)
(\sin \varphi)^{p-2} \cos \varphi > 0
\end{aligned}
\end{equation}
for all $x\in [-1,1]$.
\end{lemma}

\paragraph{\em Proof.}
Taking into consideration \eqref{est:l1},
we may use 
$$
    M =  \left(
        H - \frac{1}{(\sin\varphi)^{p - 1}} \|f\|_{L^1(-1,1)}
    \right)
$$
in \eqref{eq:Dx}. 
This ends the proof.\hfill$\blacksquare$
\par\medskip

\begin{remark}
\rm
    Assume that $f(x) < 0$ for all $x \in (-1,1)$ such that 
    \begin{equation}
    \label{ex:fneg} 
        \int_{-1}^1 f(x)\,\mathrm{d}x = - H(\sin\varphi)^{p - 1}\,.
    \end{equation}
    Then 
    \begin{multline*}
    \min\limits_{x_0\in [-1,1]}\min\limits_{x\in [x_0,1]}
    \left(
    H^{\frac{p}{p-1}}
    +
    \frac{p}{(p-1)\cos\varphi}
    \int_x^1 \Phi_{p'}\left(\int_{x_0}^{\tau} f(s)\mathrm{d}s
    \right)\mathrm{d}\tau
    \right)
    \geq
    \\
    H^{\frac{p}{p-1}} 
    -
    \frac{p}{(p-1)\cos\varphi}
    \int_{-1}^1 H^{\frac{1}{p-1}}\sin\varphi\,\mathrm{d}x
    =
    H^{\frac{p}{p-1}} 
    -
    2\frac{p}{(p-1)}
   H^{\frac{1}{p-1}}
   \tan\varphi
   =
   \\
    H^{\frac{1}{p-1}}
    \left(
        H - 2\frac{p}{(p-1)}\tan\varphi
    \right)\,.
    \end{multline*}

Thus  Hypothesis {\rm (HF)} is satisfied for any $H > 2\frac{p}{(p-1)}\tan\varphi$. Hence, we obtain the lower bound 
$u\geq -H$ by Theorem~\ref{thm:no:touch}. But the lower bound on $D(\cdot)$
obtained from this information would be dependent on the concrete weak solution $u$ of the nonlinear problem \eqref{weak_formulation}.

On the other hand, the bound \eqref{est:l1}, which is independent of $u$, is not optimal. 
Indeed, if a negative function $f$ satisfies \eqref{ex:fneg}, then
$$
   \|f\|_{L^1(-1,1)} = H(\sin\varphi)^{p - 1}\,,
$$ 
and by \eqref{est:l1} from Theorem~\ref{thm:boundness:h}, $u\geq -H$. Hence the lower bound 
of type \eqref{estimate:Dx}  independent of $u$ obtained from \eqref{eq:Dx}, 
yields that $D(x)\geq \mbox{const.} \geq 0$ for all $x\in [-1,1]$, 
but in our further analysis we need this constant to be strictly positive.
Thus it will be very interesting for practical reasons to obtain finer
estimates of type \eqref{est:l1} in order to get finer lower bounds
of type \eqref{estimate:Dx}. {\it We leave it as an interesting and technically quite complicated open problem.}
\end{remark}

\begin{remark}
{\rm
\label{rem:pl2}
Let us note that we can handle only the case $p>2$, since we are lacking an analogue of Lemma~\ref{lemma:lower:bound:int} for $1<p<2$.
We also do not know (at the time when we wrote this paper), if the spatially dependent diffusion coefficient
$D(\cdot)$ given by~\eqref{eq:D} is essentially bounded for $1<p<2$.
If not, the weak solution of the linearized problem \eqref{eq:linearization:line:prob} has to be considered in some appropriate weighted Sobolev spaces, see, e.g.,\cite{
DrabekGirgTakac2004,
DrabekGirgTakacUlm2004,  Takac2002}. {\it Obtaining similar results such as
Lemma~\ref{lemma:lower:bound:int} and
Lemma~\ref{prop:nonlinear:to:linear}
for $1<p<2$ poses
an interesting and important open problem.}
}
\end{remark}

\subsection{Additional regularity results.}
Our next goal is to find an estimate for $\|u'\|_{\infty}$.
Essential tool of this section is the method of linearization presented in previous section. For this reason, we can deal with the case $p>2$ only. We will use the fact that any weak solution $u$ to \eqref{eq_h_elliptic} satisfies $u > -H$ by Theorem~\ref{thm:no:touch} and hence for any weak solution $u$ to \eqref{eq_h_elliptic} there exist $M > 0$ such that
\begin{equation}
\label{eq:def:M}
    \min\limits_{x\in [-1,1]} \left(u(x)+H\right) \geq M > 0\,.
\end{equation}

Main result of this section is the following theorem.
\begin{theorem}
\label{thm:boundedness:derivative}
Let $p>2$, Hypothesis {\rm (HF)}  be satisfied,
and $u$ be any (bounded) weak solution to \eqref{eq_h_elliptic}. 
If $u'(1) \geq 0$, then
    $$        
        \|u'\|_{\infty} \leq \frac{1}{K\left(\varphi, H, p, \|u^-\|_{\infty}\right)}
        \left(2  +   \frac{D(1)}{H\cos\varphi(\sin \varphi)^{p-2}}\right)\|f\|_{L^{1}(-1,1)}\,,
    $$
else
    $$        
        \|u'\|_{\infty} \leq \frac{2}{K\left(\varphi, H, p, \|u^-\|_{\infty}\right)}\left(1 + \frac{D(1)}{H\cos\varphi(\sin \varphi)^{p-2}}\right)\|f\|_{L^{1}(-1,1)}\,,
    $$
where $K\left(\varphi, H, p, 
\|u^-\|_{\infty}\right)$ is the constant from Lemma~\ref{prop:nonlinear:to:linear}, inequality~\eqref{eq:lower:bound:D}.    
\end{theorem}

\begin{lemma}
\label{lem:derivatives:at:one}
Let $f \in L^1(-1,1)$ satisfy Hypothesis {\rm (HF)}
and $u$ be any (bounded) weak solution to \eqref{eq_h_elliptic}. 
If $u'(1) \geq 0$, then
$$
    u'(1) \leq  
    \frac{1}{H\left(\sin \varphi\right)^{p-2}\cos\varphi} \|f\|_{L^{1}(-1,1)}\,,
$$
else
$$
    |u'(1)| \leq \frac{2}{H \cos\varphi \left(\sin\varphi\right)^{p-2}} \|f\|_{L^{1}(-1,1)}\,.
$$
\end{lemma}

\paragraph{\em Proof.} 
We distinguish three basic cases possibly divided into subcases. 
In all cases except the first one 
we will use \eqref{eq:first:order:plus} 
evaluated at stationary point. Let us note that such point exists 
since $u\in W^{1,p}_0(-1,1)\hookrightarrow C[-1,1]$ with
$u(-1) = u(1)=0$, and $u \in C^1(-1,1)$ by Theorem~\ref{thm:C1:reg}.
Denote $S \stackrel{\rm def}{=}\{x\in(-1,1) \colon u'(x) = 0\}$. 
Then, for any $x_1 \in S$, we obtain from \eqref{eq:first:order:plus}
\begin{multline}
\label{eq:at:extrema}
H\left|u^{\prime}(1) \cos\varphi + \sin\varphi \right|^{p-2}\left(u^{\prime}(1) \cos \varphi+\sin\varphi\right)= 
\\
\left(u\left(x_{1}\right)+H\right)(\sin \varphi)^{p-1} - \int_{x_1}^1 f(\sigma)\, \mathrm{d} \sigma\,.
\end{multline}

\par\noindent{\it Case 1: $u'(1) = 0$.} Statement is satisfied trivially.
 
\par\noindent{\it Case 2: $u'(1) > 0$.} 
It follows that there exists $x_1 \in S$ satisfying $u(x_1)<0$.
Then it follows from \eqref{eq:at:extrema} that
\begin{multline*}
0 < H\left(\sin \varphi\right)^{p-2}\left(u^{\prime}(1) \cos \varphi + \sin\varphi \right) \leq
\\
H\left|u^{\prime}(1) \cos\varphi + \sin\varphi \right|^{p-2}\left(u^{\prime}(1) \cos\varphi + \sin\varphi\right) 
\\
=
\left(u\left(x_1\right)+H\right)(\sin\varphi)^{p-1} -
\int_{x_1}^1 f(\sigma)\, \mathrm{d} \sigma\,. 
\end{multline*}
Subtracting $H(\sin\varphi)^{p-1}$ from the previous inequality we get
$$
    0<
    H\left(\sin\varphi\right)^{p-2}u^{\prime}(1) \cos\varphi
    \leq 
    u\left(x_{1}\right)\left(\sin\varphi\right)^{p-1} -
    \int_{x_1}^1 f(\sigma)\, \mathrm{d} \sigma\,. 
$$
Using the fact that $u(x_1) < 0$, we finally obtain
\begin{multline}
0 < H\left(\sin\varphi\right)^{p-2} u'(1) \cos\varphi
\leq
u\left(x_{1}\right)(\sin\varphi)^{p-1} -
\int_{x_1}^1 f(\sigma) \,\mathrm{d} \sigma 
< 
\\
\int_{x_1}^1 |f(\sigma)| \,\mathrm{d} \sigma 
\leq
\|f\|_{L^{1}(-1,1)}\,. 
\end{multline}
It follows
$$
0 < u'(1) < 
\frac{1}{H\left(\sin \varphi\right)^{p-2}\cos\varphi} \|f\|_{L^{1}(-1,1)}\,.
$$

\par\noindent{\it Case 3: $u'(1) < 0$.} Now there exists $x_1 \in S$ satisfying $u(x_1) > 0$. Further, this case is divided into three subcases. 

\par\noindent{\it Subcase 3a: $0 > u'(1) \geq - \tan\varphi$.} Then $0 \leq u'(1) \cos\varphi + \sin\varphi < \sin\varphi$.
Using \eqref{eq:at:extrema} we get
\begin{multline*}
    H(\sin\varphi)^{p - 2} (u'(1) \cos\varphi + \sin\varphi)
    \\
    \geq
    H\left|u'(1)\cos\varphi + \sin\varphi\right|^{p - 2}(u'(1) \cos\varphi + \sin\varphi)
    \\
    =
    (u(x_{1}) + H)(\sin\varphi)^{p - 1} - \int_{x_{1}}^1 f(\sigma)\,\mathrm{d}\sigma\,.
\end{multline*}
Hence
$$
    0 > H\cos\varphi(\sin\varphi)^{p - 2} u'(1) \geq u(x_1)(\sin\varphi)^{p - 1} 
    - 
    \int_{x_1}^1 f(\sigma)\,\mathrm{d}\sigma > - 
    \int_{x_1}^1 f(\sigma)\,\mathrm{d}\sigma \,.
$$
It follows
$$
|u'(1)| < \frac{1}{H\left(\sin \varphi\right)^{p-2}\cos\varphi} \|f\|_{L^{1}(-1,1)}\,.
$$
\par\noindent{\it Subcase 3b: $-\tan\varphi > u'(1) > - 2\tan\varphi$.} 
Then 
$$
    0 > u'(1) \cos\varphi + \sin\varphi > -\sin\varphi
$$
and 
$$
    (\sin\varphi)^{p - 2} > \left|u(1)\cos\varphi + \sin\varphi\right|^{p - 2}\,.
$$
Using \eqref{eq:at:extrema} 
where we subtracted $H (\sin\varphi)^{p - 1}$ from both sides 
and equality~\eqref{eq:odd:powers:equality}, we get
\begin{multline*}
    0 > H\left(
    \left|u'(1)\cos\varphi + \sin\varphi\right|^{p - 2}(u'(1) \cos\varphi + \sin\varphi) 
    - \left|\sin\varphi\right|^{p - 2}\sin\varphi 
    \right)
    \\
    =
    H(p - 1)\int_0^1|\sin\varphi + \theta(u'(1)\cos\varphi + \sin\varphi)|^{p - 2}\,\mathrm{d}\theta\,
    u'(1)\cos\varphi
    \\
    =
    u(x_1)(\sin\varphi)^{p - 1} - \int_{x_1}^1 f(\sigma)\,\mathrm{d}\sigma\,.
\end{multline*}
Hence
\begin{multline*}
    0 <  H(p - 1)\int_0^1|\sin\varphi + \theta(u'(1)\cos\varphi + \sin\varphi)|^{p - 2}\,\mathrm{d}\theta\,
   (-u'(1))\cos\varphi
    \\
    =
    -u(x_1)(\sin\varphi)^{p - 1} + \int_{x_1}^1 f(\sigma)\,\mathrm{d}\sigma <  \int_{x_1}^1 f(\sigma)\,\mathrm{d}\sigma\,,
\end{multline*}
which follows
$$
    0 < |u'(1)| \leq \frac{
    \|f\|_{L^1(-1,1)}
}
    {H\cos\varphi(p - 1) \int_0^1|\sin\varphi + \theta(u'(-1)\cos\varphi + \sin\varphi)|^{p - 2}\,\mathrm{d}\theta}\,.
$$
Since  
$$
    (p - 1)\int_0^1|1 + \theta(u'(1)\cot\varphi + 1)|^{p - 2}\,\mathrm{d}\theta
    \geq
   \frac{1}{2}
$$
by Lemma~\ref{lemma:lower:bound:int}, we obtain an upper bound
$$
    0 < |u'(1)| 
    \leq 
    2\frac{\|f\|_{L^1(-1,1)}}{H\cos\varphi (\sin\varphi)^{p - 2}}\,.
$$
%The last inequality follows from \cite[Lemma~A.1, p.~233]{Takac2002} (see \eqref{eq:Dx} for more details).

\par\noindent{\it Subcase 3c: $-2\tan\varphi \geq u'(1) $.} 
Then 
$$ 
    -\sin\varphi \geq u'(1) \cos\varphi + \sin\varphi
$$
and 
$$
(\sin\varphi)^{p - 2} \leq \left|u(1)\cos\varphi + \sin\varphi\right|^{p - 2}\,.
$$
Using \eqref{eq:at:extrema} we get
\begin{multline*}
    0 > H(\sin\varphi)^{p - 2} (u'(1) \cos\varphi + \sin\varphi)
    \\
    \geq
    H\left| u'(1) \cos\varphi + \sin\varphi \right|^{p - 2}  (u'(1) \cos\varphi + \sin\varphi)
    \\
    =
    (u(x_1) + H)(\sin\varphi)^{p - 1} - \int_{x_1}^{1} f(\sigma)\,\mathrm{d}\sigma\,.
\end{multline*}
The rest of the proof is similar to the one of Subcase 3a.
\par\hfill$\blacksquare$
\par\medskip

\paragraph{\em Proof of Theorem \ref{thm:boundedness:derivative}.}
We will use linearization at the trivial solution described in Section~\ref{par:linearizaton}. Thus
an arbitrary but fixed weak solution $u$ to \eqref{eq_h_elliptic} also satisfies the following weak formulation of 
the following linear problem
\begin{equation}
\label{eq:linear:in:-1}
\begin{aligned}
&\int_{-1}^{1} \left( D(x) u^{\prime}(x)
+(\sin \varphi)^{p-1} u(x) - \int_x^1 f(\sigma) \mathrm{d} \sigma\right) v^{\prime}(x)\, \mathrm{d} x= 0
\end{aligned}
\end{equation}
with the diffusion coefficient given by \eqref{eq:D} for the fixed weak solution $u$ to \eqref{eq_h_elliptic}.
Hence
\begin{equation}
\label{eq:linear:with:kappa}
\begin{aligned}
& D(x) u^{\prime}(x)
+(\sin \varphi)^{p-1} u(x) - \int_x^1 f(\sigma) \,\mathrm{d} \sigma = \kappa\,,
\end{aligned}
\end{equation}
where constant
$$
\begin{aligned}
& \kappa = D(1) u^{\prime}(1)
+(\sin \varphi)^{p-1} u(1) +\cos\varphi \int_{1}^{1} f(\sigma)\, \mathrm{d} \sigma = D(1) u^{\prime}(1)
\end{aligned}
$$
by taking $x = 1$ in the previous equation.
We substitute $D(1) u^{\prime}(1)$ for $\kappa$ in \eqref{eq:linear:with:kappa} to obtain 
$$
\begin{aligned}
& D(x) u^{\prime}(x) =
-(\sin \varphi)^{p-1} u(x) + \int_x^1 f(\sigma)\, \mathrm{d} \sigma +  D(1) u^{\prime}(1)\,.
\end{aligned}
$$
Thus
\begin{equation}
\label{up:common:est}
\begin{aligned}
& |u^{\prime}(x)| \leq
\frac{1}{\min\limits_{x\in [-1,1]}D(x)}\left(
\sin (\varphi)^{p-1} |u(x)| 
+  \|f\|_{L^1(-1,1)}
+  D(1) |u^{\prime}(1)|
\right)\,.
\end{aligned}
\end{equation}
by triangle inequality.
For $u'(1)\geq 0$, we have
\begin{multline}
\label{eq:common:est:geq0}
|u^{\prime}(x)| \leq \\
\frac{1}{\min\limits_{x\in [-1,1]}D(x)}\left(
(\sin \varphi)^{p-1} \frac{1}{(\sin \varphi)^{p-1}} 
+ 1 
+   \frac{D(1)}{H\cos\varphi(\sin \varphi)^{p-2}}\right)\|f\|_{L^{1}(-1,1)}
\end{multline}
by Theorem~\ref{thm:boundness:h} and Lemma~\ref{lem:derivatives:at:one}.
Using \eqref{eq:lower:bound:D} we finally obtain
$$
\begin{aligned}
& |u^{\prime}(x)| \leq
\frac{1}{K\left(\varphi, H, p, \|u^-\|_{\infty}\right)}\left(
2 
+   \frac{D(1)}{H\cos\varphi(\sin \varphi)^{p-2}}\right)\|f\|_{L^{1}(-1,1)}\,.
\end{aligned}
$$
For $u'(1)<0$, we have different bound on $|u'(1)|$ and hence we obtain
\begin{multline}
\label{eq:common:est:l0}
|u^{\prime}(x)| \leq \\
\frac{1}{\min\limits_{x\in [-1,1]}D(x)}\left(
(\sin \varphi)^{p-1} \frac{1}{(\sin \varphi)^{p-1}} 
+ 1 
+   \frac{2 D(1)}{H\cos\varphi(\sin \varphi)^{p-2}}\right)\|f\|_{L^{1}(-1,1)}
\end{multline}
by Theorem~\ref{thm:boundness:h} and Lemma~\ref{lem:derivatives:at:one}.
Then we have
$$
\begin{aligned}
& |u^{\prime}(x)| \leq
\frac{2}{K\left(\varphi, H, p, \|u^-\|_{\infty}\right)}\left(
1  
+   \frac{ D(1)}{H\cos\varphi(\sin \varphi)^{p-2}}\right)\|f\|_{L^{1}(-1,1)}
\end{aligned}
$$
using \eqref{eq:lower:bound:D} again.
\hfill$\blacksquare$
\par\medskip

The estimate on $\|u'\|_{\infty}$ from the previous theorem can be further refined to become a~priori bound (independent of $\|u^-\|_{\infty}$ and of $D(1)$) under an additional condition.

\begin{theorem}
\label{thm:apriory:bound:derivative}
Let $p>2$, Hypothesis {\rm (HF)}  be satisfied, and $u$ be any (bounded) weak solution to \eqref{eq_h_elliptic}.
Moreover, let there exists $\beta>0$
such that
\begin{equation}
\label{pom:f}
   \|f\|_{L^1(-1,1)} \leq \beta < 
   H(\sin\varphi)^{p - 1}\,.
\end{equation}
Then there exists 
a~constant $C > 0$ depending only on  $\varphi, H, p$, and $\beta$ such that 
\begin{equation}
    \label{eq:u:prime:uniform:bound}
        \|u'\|_{\infty} \leq C(\varphi, H, p, \beta) \|f\|_{L^{1}(-1,1)}\,.
\end{equation}  
\end{theorem}
\paragraph{\em Proof.} 
We proceed as in the proof of Theorem~\ref{thm:boundedness:derivative} and arrive at \eqref{up:common:est}. Then for $u'(1)\geq 0$, as in the proof of Theorem~\ref{thm:boundedness:derivative} we obtain \eqref{eq:common:est:geq0}, that is,
\begin{equation}
\label{pom:unif:d:1}
|u^{\prime}(x)| \leq 
\frac{1}{\min\limits_{x\in [-1,1]}D(x)}\left(
2 
+   \frac{D(1)}{H\cos\varphi(\sin \varphi)^{p-2}}\right)\|f\|_{L^{1}(-1,1)}\,.
\end{equation}
For $u'(1)<0$, as in the proof of Theorem~\ref{thm:boundedness:derivative} we obtain \eqref{eq:common:est:l0}, that is,
\begin{equation}
\label{pom:unif:d:2}
|u^{\prime}(x)| \leq
\frac{2}{\min\limits_{x\in [-1,1]}D(x)}\left(
1 
+ 
\frac{D(1)}{H\cos\varphi(\sin \varphi)^{p-2}}\right)\|f\|_{L^{1}(-1,1)}\,.
\end{equation}

Under assumption \eqref{pom:f}, we have
\begin{equation}
\label{pom:d}
D(x) \geq K'\left( \varphi, H, p, \|f\|_{L^1(-1,1)}\right)>0
\end{equation}
by estimate~\eqref{estimate:Dx} from Lemma~\ref{prop:uniform:lower:bound}. Using \eqref{pom:unif:d:1}, \eqref{pom:unif:d:2} and \eqref{pom:d}, we find
\begin{equation}
\label{pom:unif:d:3}
|u^{\prime}(x)| \leq
\frac{2}{K'\left( \varphi, H, p, \|f\|_{L^1(-1,1)}\right)}\left(
1 
+ 
\frac{D(1)}{H\cos\varphi(\sin \varphi)^{p-2}}\right)\|f\|_{L^{1}(-1,1)}\,.
\end{equation}

Now, using boundary condition $u(1) = 0$, triangle inequality,  $x^s < x$ provided $x > 1$ for $0 < s < 1$ and
$(a + b)^s \leq 2^{s-1}(a^s + b^s)$ provided $a$,$b \geq 0$ for $s \geq 1$,
and Lemma~\ref{lem:derivatives:at:one}, we obtain
%\begin{multline}
\begin{equation}
\label{eq:D:1:bound}
\begin{aligned}
0 < D(1) 
     & 
    = (u(1)+H)(p-1) \int_{0}^{1}\left|\sin \varphi+\theta u^{\prime}(1) \cos \varphi\right|^{p-2} \mathrm{~d} \theta \cos \varphi 
    \\   & 
    \leq (p - 1) H \max\{1, 2^{p - 1}\}\int_{0}^{1} \left(|\sin \varphi|^{p - 2} + \left|\theta u^{\prime}(1) \cos \varphi \right|^{p-2} \right)\mathrm{~d} \theta \cos \varphi \\    & 
    = (p - 1) H  \max\{1, 2^{p - 1}\}\left((\sin \varphi)^{p - 2} + \frac{\left|u^{\prime}(1) \cos \varphi \right|^{p-2}}{p - 1} \right)\cos\varphi
     \\    &
    \leq (p - 1) H \max\{1, 2^{p - 1}\} \cos\varphi (\sin \varphi)^{p - 2} 
    \\
    &
    + 2\max\{1, 2^{p - 1}\}\frac{(\cos\varphi)^{p - 2}}{(\sin\varphi)^{p - 2}} \|f\|_{L^{1}(-1,1)}\,. 
\end{aligned}
\end{equation}
Setting $$C'(\varphi, H, p, \beta)
\stackrel{\mathrm{def}}{=} (p - 1) H \max\{1, 2^{p - 1}\} \cos\varphi (\sin \varphi)^{p - 2} + 2\max\{1, 2^{p - 1}\}\frac{(\cos\varphi)^{p - 2}}{(\sin\varphi)^{p - 2}} \beta$$
and taking into account that
$$
0< 
K'\left( \varphi, H, p, \beta\right)
\leq
K'\left( \varphi, H, p, \|f\|_{L^1(-1,1)}\right)\,,
$$
we deduce from \eqref{pom:unif:d:3} and \eqref{eq:D:1:bound} that
\begin{equation}
\label{pom:unif:d:4}
|u^{\prime}(x)| \leq
\frac{2}{K'\left( \varphi, H, p, \beta\right)}\left(
1 
+ 
\frac{C'(\varphi, H, p, \beta)}{H\cos\varphi(\sin \varphi)^{p-2}}\right)\|f\|_{L^{1}(-1,1)}\,.
\end{equation}
This establishes a~priori bound \eqref{eq:u:prime:uniform:bound}.
\hfill$\blacksquare$
\par\medskip

\subsection{Weak and Strong Maximum Principles via linearization}

In this section, we derive Weak and Strong Maximum Principles for \eqref{eq_h_elliptic}. 
\begin{comment}
\color{red}
We will make use of the following standard notation here
$$
u^+\stackrel{\mathrm{def}}{=}\max\{u,0\}\,,\quad
u^-\stackrel{\mathrm{def}}{=}\max\{-u,0\}
$$
for 
$u\in W_0^{1,p}(-1,1)$.
Then $u = u^+ - u^-$ and let us recall that also
$u^-, u^+\in W_0^{1,p}(-1,1)$
(see, e.g.,
\cite[Coroll. 2.1.8, p.~47]{Ziemer89}).
Nyni bych psal:
\end{comment}
Let us recall that,
for 
$u\in W_0^{1,p}(-1,1)$,
functions
$$
u^+\stackrel{\mathrm{def}}{=}\max\{u,0\}\,,\quad
u^-\stackrel{\mathrm{def}}{=}\max\{-u,0\}
$$
also satisfy
$u^-, u^+\in W_0^{1,p}(-1,1)$
(see, e.g.,
\cite[Coroll. 2.1.8, p.~47]{Ziemer89}).

It is inherent to assume that $f \geq 0$ a.e. when studying Weak or Strong Maximum Principle. 
Let us point out that the assumption $f \geq 0$ a.e. implies that 
Hypothesis {\rm (HF)} is satisfied. 

Note that even Weak Maximum Principle for the problem \eqref{eq_h_elliptic} is not as straightforward to prove as in the case of the classical $p$-Laplacian problem, where one can use $u^-$ as a test function in the weak formulation to prove the result. Here, the problem is caused by the term $u(x)+H$, which multiplies the part depending on the derivative of the solution. This causes serious difficulty in finding a suitable choice of test function, which would lead to conclusion. Fortunately, the linearization method from Section~\ref{par:linearizaton} can be used under certain rather general assumptions.
This procedure yields the following statement.

\begin{proposition}
[Weak Maximum Principle]
\label{wmp}
Let $p>2$, $f \in L^1(-1,1)$, and
$u \in W_0^{1,p}(-1,1)$ be a weak solution to \eqref{eq_h_elliptic}. 
If $f \geq 0$ a.e. in $(-1,1)$, then $u \geq 0$ in $(-1,1)$.
\end{proposition}
\paragraph{\em Proof.}
Let $u\in W_0^{1,p}(-1,1)$
be any weak solution to \eqref{eq_h_elliptic}.
Then once this solution is fixed, we define function $D(\cdot)$
by
\eqref{eq:D}. Now, by  Lemma~\ref{prop:nonlinear:to:linear},  this same solution $u$ also satisfies
\eqref{eq:weak:linear}, that is,
$$
\int_{-1}^{1} D(x) u'(x) v'(x)\,\mathrm{d}x + 
(\sin \varphi)^{p - 1} \int_{-1}^{1} u(x) v'(x) \,\mathrm{d}x 
=
\int_{-1}^{1} f(x) v(x) \,\mathrm{d}x
$$
for any $v \in W^{1,p}_0(-1,1)$.
Taking $v = u^-$ we have
\begin{multline}
\label{eq:wmp:testing}
    -\int_{-1}^{1} D(x)\left( \left[u^-(x)\right]'\right)^2 \,\mathrm{d}x
    - 
    (\sin \varphi)^{p - 1} \int_{-1}^{1} 
    \frac{1}{2}\left[\left(u^-(x)\right)^2\right]'  \,\mathrm{d}x 
    \\
    =
    \int_{-1}^{1} f\ [u(x)]^- \,\mathrm{d}x
    \geq 0\,.
\end{multline}
Since
$$D(x)\geq
K\left(\varphi, H, p, \|u^-\|_{\infty}\right)
>0
$$
for all $x\in [-1,1]$
by {
Lemma~\ref{prop:nonlinear:to:linear} and
$$
\int_{-1}^{1} 
    \frac{1}{2}\left[\left(u^-(x)\right)^2\right]'  \,\mathrm{d}x
    =
    \frac{1}{2}
    \left[
    \left(u^-(1)\right)^2-
    \left(u^-(-1)\right)^2
    \right]
    =0
$$
by boundary conditions
$u(-1)=0=u(1)$,
we infer from \eqref{eq:wmp:testing}
that
$$
K\left(\varphi, H, p, \|u^-\|_{\infty}\right)
\int_{-1}^{1}\left( \left[u^-(x)\right]'\right)^2 \,\mathrm{d}x\leq
\int_{-1}^{1} D(x)\left( \left[u^-(x)\right]'\right)^2 \,\mathrm{d}x\leq 0\,.
$$
Hence $[u^-]' = 0$ in $(-1,1)$ and consequently
$u^- = 0$ in $(-1,1)$ since $u^-(-1) = 0 = u^-(1)$ (see, e.g.,
\cite[Coroll. 2.1.9, p.~47]{Ziemer89}).
Thus $u=u^+\geq 0$ in $(-1,1)$.
\hfill$\blacksquare$
\par\medskip

As a useful consequence of Weak Maximum Principle, we obtain the following result.
\begin{corollary}
\label{cor:D:positive}
Let $p>2$, $f \in L^1(-1,1)$,
and $f\geq 0$.
Let
$u \in W_0^{1,p}(-1,1)$ be a weak solution to \eqref{eq_h_elliptic}.  Then 
\begin{equation}
\label{estimate:Dx:pos}
D(x) \geq  
\frac{1}{2}
H (\sin \varphi)^{p-2} \cos \varphi > 0
\end{equation}
for all $x\in [-1,1]$.
\end{corollary}

\paragraph{\em Proof.}
By Proposition \ref{wmp}, $u\geq 0$.
Hence,
we can use $M=H$
in 
\eqref{eq:Dx}.
This establishes 
\eqref{estimate:Dx:pos}.
\hfill$\blacksquare$
\par\medskip

Since $D(\cdot)\in C[-1,1]$ by Lemma~\ref{prop:nonlinear:to:linear}, we get from~\eqref{estimate:Dx:pos} that the function $1/D(\cdot)\in C[-1,1]$ and hence it is integrable. 

\par\medskip

The crucial part of our proof of the validity of Strong Maximum Principle 
is to show that Green's function associated with the operator
$$
\mathop{\mathcal{L}}u(x) \stackrel{\rm def}{=} - \left(D(x) u'(x) \right)' 
 - (\sin \varphi)^{p - 1} u'(x)
$$
is positive. 
At first, we prove that
\begin{equation}
\label{eq:def:greenka}
G(x,x_0) \stackrel{\rm def}{=} \left\{
\begin{aligned}
    &\frac{E^+(x_0) - 1}{E^-(x_0) - E^+(x_0)} \frac{1}{(\sin \varphi)^{p - 1}} 
    \left(
        E^-(x) 
        - 1
    \right)
    && \text{ for } x \in [-1, x_0]\,,\\
    & \frac{E^-(x_0) - 1}{E^-(x_0) - E^+(x_0)} \frac{1}{(\sin \varphi)^{p - 1}} 
    \left(
        E^+(x)
        - 1
    \right)
    && \text{ for } x \in (x_0, 1]\,,
\end{aligned}
\right.
\end{equation}
$x_0 \in [-1,1]$, is indeed the Green's function associated with $\mathop{\mathcal{L}}$.
Here,
$$
    E^-(s) \stackrel{\rm def}{=}  \exp\left[-\int_{-1}^{s} \frac{(\sin \varphi)^{p - 1}}{D(\xi)}\,\mathrm{d}\xi\right]\,, 
$$
and
$$
    E^+(s) \stackrel{\rm def}{=} \exp\left[\int_{s}^{1} \frac{(\sin \varphi)^{p - 1}}{D(\xi)}\,\mathrm{d}\xi\right]\,.
$$
Note that both $E^-(s)$ and $E^+(s)$ make sense since $D(\xi)\geq \mathrm{const.} > 0$ for $\xi\in[-1,1]$ by Corollary~\ref{cor:D:positive} and it is continuous by Lemma~\ref{prop:nonlinear:to:linear}.} 
Moreover,
\begin{multline}
\label{EpminEm}
E^+(s)-E^-(s) = \\
\exp\left[\int_{s}^{1} \frac{(\sin \varphi)^{p - 1}}{D(\xi)}\,\mathrm{d}\xi\right]
\left(1-\exp\left[-\int_{-1}^{1} \frac{(\sin \varphi)^{p - 1}}{D(\xi)}\,\mathrm{d}\xi\right]\right)\geq \\ \geq
\left(1-\exp\left[-\int_{-1}^{1} \frac{(\sin \varphi)^{p - 1}}{D(\xi)}\,\mathrm{d}\xi\right]\right) = 
\mathrm{const.}>0\,.
\end{multline}

\begin{comment}
With this at hand, 
it is a matter of straightforward calculation to show that for all $x_0 \in (-1,1)$ we have:
\begin{itemize}
    \item[(G1)] $G(\cdot,x_0) \in C[-1,1]$,
    \item[(G2)] $G(\cdot,x_0) \in C^{\infty}(-1,x_0)$,
    \item[(G3)] $G(\cdot,x_0) \in C^{\infty}(x_0, 1)$,
    \item[(G4)] $G_x(\cdot,x_0) \in L^{\infty}(-1, 1)$, and
    \item[(G5)] $G_{xx}(\cdot,x_0) \in L^{\infty}(-1, 1)$.
\end{itemize}
Recall that $D(x) \geq K > 0$.
\end{comment}

The following result is very usefull in establishing that
$G$ is Green's function associated with $\mathop{\mathcal{L}}$.
\begin{lemma}
\label{lem:GL}
The function $G\colon [-1,1]\times[-1,1]\to\mathbb{R}$ given by \eqref{eq:def:greenka}
is continuous and satisfies that
\begin{itemize} 
\item[(GL)]
there exists
$\kappa>0$ such that $|G(s,y)-G(t,y)|\leq \kappa |s-t|$ for any $s,t,y\in [-1,1]$.
\end{itemize}
Moreover, the partial derivative $G_x(x,y)$
exists for all $(x,y)\in (-1,1)\times(-1,1)$ such that
$x\not=y$.
\end{lemma}

\paragraph{\em Proof.}
It is easy to see from \eqref{eq:def:greenka} that
$G\in C([-1,1]\times[-1,1])$ and that the partial derivative $G_x=\partial G/\partial x$ exists for all $(x,y)\in (-1,1)\times (-1,1)$ such that $x\not=y$. 
Then it is straightforward to verify that
$G\in C^1(\overline{\Omega_1})$ and $G\in C^1(\overline{\Omega_2})$, where
\begin{eqnarray}
\Omega_1&\stackrel{\mathrm{def}}{=}&
\{(x,y)\in\mathbb{R}\colon -1<x<1\mbox{ and } 0<y<x\}\,, \\
\Omega_2&\stackrel{\mathrm{def}}{=}&
\{(x,y)\in\mathbb{R}\colon -1<x<1\mbox{ and } x<y<1\}\,.
\end{eqnarray}
Thus, for all $y\in [-1,1]$ fixed, the function $G(\cdot, y)$ is absolutely continuous
 (since continuous gluing of two absolutely continuous functions is absolutely continuous).
Now, let us set 
\begin{equation}
\label{def:emko}
\kappa= 
\max_{i=1,2}\left(
\max_{(x,y)\in
\overline{\Omega_i}}\left|\vphantom{\sum}G_x(x,y)\right|\right)\,.
\end{equation}
Then $|G_x(x,y)|\leq \kappa$
for all 
$x,y\in [-1,1]$ such that $x\not=y$. Now let us consider $y\in [-1,1]$ arbitrary but fixed. Absolute continuity of $G(\cdot, y)$ together with the fact that $|G_x(\cdot, y)|<\kappa$ a.e. in $[-1,1]$ imply that $G(\cdot, y)$ is Lischitz continuous with constant $\kappa$ defined by \eqref{def:emko}. As $\kappa$ is independent of $y\in [-1,1]$,
we established the  statement of the lemma.
\hfill$\blacksquare$
\par\medskip

\begin{theorem}[Green's function]
\label{thm:conv}
Let $f$ satisfy {\rm (HF)}.
The function 
    $$
        u(x) = \int_{-1}^1 G(x,y) f(y)\,\mathrm{d}y
    $$
    is weak solution to the linear problem \eqref{eq:linearization:line:prob}.
\end{theorem}
\paragraph{\em Proof.}
Since 
$$
    G(-1,y) = \frac{E^+(y) - 1}{E^-(y) - E^+(y)} \frac{1}{(\sin \varphi)^{p - 1}} 
    \left(
        \mathrm{e}^0 
        - 1
    \right)
    = 0
$$
and
$$
    G(1,y) = \frac{E^-(y) - 1}{E^-(y) - E^+(y)} \frac{1}{(\sin \varphi)^{p - 1}} 
    \left(
        \mathrm{e}^0
        - 1
    \right)
    = 0
$$
for all $y \in (-1,1)$, we have
$$
    u(-1) = \int_{-1}^1 G(-1,y) f(y)\,\mathrm{d}y = 0
$$
and
$$
     u(1) = \int_{-1}^1 G(1,y) f(y)\,\mathrm{d}y = 0\,.
$$
Thus $u$ satisfies boundary conditions.

Let $x\in (-1,1)$ be arbitrary. Then
\begin{align}
\nonumber
&
\frac{\mathrm{d}u}{\mathrm{d}x}(x)= 
\frac{\mathrm{d}}{\mathrm{d}x} \int_{-1}^1 G(x,y) f(y)\,\mathrm{d}y 
=
\\ &
\nonumber
=
\lim\limits_{\triangle x\to 0}
\frac{
\int_{-1}^1 G(x+\triangle x,y) f(y)\,\mathrm{d}y-
\int_{-1}^1 G(x,y) f(y)\,\mathrm{d}y}{\triangle x} 
=
\\ &
\label{eq:b1}
=
\lim\limits_{\triangle x\to 0}
\int_{-1}^1
\frac{(G(x+\triangle x,y)-G(x,y) f(y)}{\triangle x}\,\mathrm{d}y
%\\ &
%=
%\lim\limits_{\triangle x\to %0}\,
%\int_{-1}^{x-\varepsilon}\,+\,\int_{x-\varepsilon}^{x+\varepsilon}\,+\,\int_{x+\varepsilon}^1\quad 
%\frac{(G(x+\triangle x,y)-G(x,y) f(y)}{\triangle x}\,\mathrm{d}y\,.
\end{align}
Let us note that, for any fixed $x\in (-1,1)$, the partial derivative $$G_x(x,y)=\lim\limits_{\triangle x\to 0}\frac{G(x+\triangle x,y)-G(x,y)}{\triangle x}$$ is well defined for any  $y\in (-1,1)\setminus\{x\}$ by Lemma~\ref{lem:GL}.
%At first, let us note that $G_x(x,y)$
%exists on $[x, x-\varepsilon]\times [-1,1]$.
Now taking into account (GL), we can use function $\kappa |f|\in L^1(-1,1)$ as an integrable majorant for the integrand in \eqref{eq:b1}. Then, by the Lebesgue dominated convergence theorem,  
we obtain
\begin{multline}
\label{eq:b2}
\lim\limits_{\triangle x\to 0}
\int_{-1}^{1}
\frac{
(G(x+\triangle x,y)-G(x,y) f(y)}{\triangle x}\,\mathrm{d}y
=\\=
\int_{-1}^{1}
\lim\limits_{\triangle x\to 0}
\frac{
(G(x+\triangle x,y)-G(x,y) f(y)}{\triangle x}\,\mathrm{d}y=
\int_{-1}^{1}
G_x(x,y) f(y)\,\mathrm{d}y\,.
\end{multline}
Thus, we obtain
\begin{equation}
\label{eq:b5}
\frac{\mathrm{d}u}{\mathrm{d}x}(x)=
\int_{-1}^{1}
G_x(x,y) f(y)\,\mathrm{d}y
\end{equation}
for any $x\in (-1,1)$.
Now, using formulas for $G$ from \eqref{eq:def:greenka} and taking them into derivative, we obtain
\begin{multline}
\label{eq:dudx:B}
    \frac{\mathrm{d}u}{\mathrm{d}x}(x) 
    = \\ =
\frac{E^+(x)}{D(x)}
\int_{-1}^x \frac{E^-(y) - 1}{E^+(y) - E^-(y)}
         f(y) \,
    \mathrm{d}y
    +\frac{E^-(x)}{D(x)}
    \int_x^1 \frac{E^+(y) - 1}{E^+(y) - E^-(y)}
         f(y) \,
    \mathrm{d}y\,,
\end{multline}
where we used
$$    \frac{\mathrm{d}E^{\pm}}{\mathrm{d}s}(s) = - \frac{(\sin\varphi)^{p - 1}}{D(s)} E^{\pm}(s)\,.
$$
Let us note that the function on the right-hand side of \eqref{eq:dudx:B} can be continuously extended to $[-1,1]$ as its limits exist and are finite on both ends of the interval. Thus $u\in C^1[-1,1]$ and hence
its weak and strong derivatives coincide
and $u\in W_0^{1,p}(-1,1)$.

From \eqref{eq:dudx:B}, we obtain
\begin{equation}
\label{eq:dudx:mod}
\begin{aligned}
    &
D(x)\frac{\mathrm{d}u}
    {\mathrm{d}x}(x) 
    = 
    E^+(x)\int_{-1}^x \frac{E^-(y) - 1}{E^+(y) - E^-(y)} f(y) \,
    \mathrm{d}y
    \\
    &
    +
    E^-(x) \int_x^1 \frac{E^+(y) - 1}{E^+(y) - E^-(y)} 
    f(y) \,
    \mathrm{d}y
\end{aligned}
\end{equation}
for all $x\in (-1,1)$.
Taking into account
\eqref{EpminEm}, and definitions of $E^{\pm}$, and the fact that $f\in L^1(-1,1)$, we deduce that the integrands on the right-hand side of \eqref{eq:dudx:mod} are from $L^1(-1,1)$. Thus the continuous extension of $D(\cdot)\,\mathrm{d}u/{\mathrm{d}x}$ to $[-1,1]$ belongs to $AC[-1,1]$. Hence the derivative of $D(\cdot)\,\mathrm{d}u/{\mathrm{d}x}$
exists a.e. in~$(-1,1)$ and can be obtained by a straightforward calculation (as it contains products and sums of functions from $AC[-1,1]$). 
Indeed, 
we obtain
\color{black}
\begin{equation}    \label{eq:verification:linearized:laplace}
\begin{aligned}
    &
    -
    \frac{\mathrm{d}}{\mathrm{d}x}\left(D(x)\frac{\mathrm{d} u}{\mathrm{d}x}(x)\right)
    =
    \frac{E^-(x) -1}{E^-(x) - E^+(x)}E^+(x) f(x)
    \\
    &
    -
    \frac{E^+(x) -1}{E^-(x) - E^+(x)}E^-(x) f(x)
    \\
    &
    +
    \frac{(\sin \varphi)^{p - 1}E^+(x)}{D(x)}
    \int_{-1}^x 
        \frac{E^-(y) -1}{E^+(y) - E^-(y)} f(y)\,
    \mathrm{d}y
    \\
    &
    +
   \frac{(\sin \varphi)^{p - 1} E^-(x)}{D(x)}
    \int_x^1
        \frac{E^+(y) -1}{E^+(y) - E^-(y)}
        f(y)\,
    \mathrm{d}y
\end{aligned}
\end{equation}
a.e. in $(-1,1)$.
Using \eqref{eq:dudx:B}, we obtain
\begin{equation}
\label{eq:verification:linearized:drift}
\begin{aligned}
    &
    (\sin\varphi)^{p - 1} \frac{\mathrm{d} u}{\mathrm{d}x}(x)
    =
\frac{(\sin \varphi)^{p - 1}E^+(x)}{D(x)}    
\int_{-1}^x 
        \frac{E^-(y) -1}{E^+(y) - E^-(y)} f(y)\,
    \mathrm{d}y
    \\
    &
    +
    \frac{(\sin \varphi)^{p - 1}E^-(x)}{D(x)}
    \int_x^1
        \frac{E^+(y) -1}{E^+(y) - E^-(y)}
        f(y)\,
    \mathrm{d}y\,.
\end{aligned}
\end{equation}
Subtracting \eqref{eq:verification:linearized:drift} from \eqref{eq:verification:linearized:laplace},
we see that
\begin{align*}
    &
    - 
    \frac{\mathrm{d}}{\mathrm{d}x}\left(D(x)\frac{\mathrm{d} u}{\mathrm{d}x}(x)\right)
    -  
    (\sin\varphi)^{p - 1} \frac{\mathrm{d} u}{\mathrm{d}x}(x)
    \\
    &
    =
    \frac{E^-(x) -1}{E^-(x) - E^+(x)}E^+(x) f(x)
    -
    \frac{E^+(x) -1}{E^-(x) - E^+(x)}E^-(x) f(x)
    \\
    &
    =
    \frac{E^-(x)E^+(x) - E^+(x) - E^-(x)E^+(x) + E^-(x)}{E^-(x) - E^+(x)} f(x)
    =
    f(x)
\end{align*}
for almost every $x \in (-1,1)$.
Now, let us recall that
the continuous extension of $D(\cdot)\,\mathrm{d}u/{\mathrm{d}x}$ to $[-1,1]$ belongs to $AC[-1,1]$. Thus, for any 
$v\in W^{1,p}_0(-1,1)$
(working with its representative in    $AC[-1,1]$), the following integration by parts makes sense
\begin{multline}
\int_{-1}^1 
D(x)u'(x)v'(x)
\mathrm{d}x
+
(\sin\varphi)^{p-1}\int_{-1}^1
u(x)v'(x)\mathrm{d}x
= \\ =
\int_{-1}^1 
\left((D(x)u'(x))'v(x)
+
(\sin\varphi)^{p-1}
u(x)'\right)v(x)\mathrm{d}x
=
\int_{-1}^1 
f(x)v(x)\mathrm{d}x\,,
\end{multline}
where we use shorter notation for derivatives.
This establishes that the function $
        u(x) = \int_{-1}^1 G(x,y) f(y)\,\mathrm{d}y
$
is a weak solution to the linear problem \eqref{eq:linearization:line:prob}.
\hfill$\blacksquare$
\par\medskip

\begin{lemma}[Positivity of Green's funtion]
\label{lem:positivity}
The Green's function
of 
$
\mathop{\mathcal{L}}
$
%is positive.
satisfies $G(x,y)>0$ for all $x,y\in (-1,1)$.
\end{lemma}
\paragraph{\em Proof.} Taking into account that 
$D(\xi)\geq\mathrm{const.}>0$, we have
\begin{eqnarray*}
0 < E^-(s) =  \exp\left[-\int_{-1}^{s} \frac{(\sin \varphi)^{p - 1}}{D(\xi)}\,\mathrm{d}\xi\right]
< 1 \\  
1 < E^+(s) = \exp\left[\int_{s}^{1} \frac{(\sin \varphi)^{p - 1}}{D(\xi)}\,\mathrm{d}\xi\right]
\end{eqnarray*}
for all $s\in (-1,1)$.
Then the positivity of Green's function given by~\eqref{eq:def:greenka} follows  
from the inequalities above combined with inequality \eqref{EpminEm}.
\hfill$\blacksquare$
\par\medskip

\begin{theorem}[Strong Maximum Principle]
\label{smp}
Let $u \in W_0^{1,p}(-1,1)$ be a solution to 
\eqref{eq_h_elliptic} with $f \in L^1(-1,1)$.
If $f \geq 0$ and $f \not\equiv 0$, then $u > 0$ on $(-1,1)$.
\end{theorem}
\paragraph{\em Proof.}
We showed in Section~\ref{par:linearizaton}, 
that any solution $u$ to nonlinear problem \eqref{eq_h_elliptic}
is also solution to linear problem \eqref{eq:linearization:line:prob} 
with spatially dependent diffusion coefficient $D(\cdot)$ constructed from $u$. 
By our assumption 
there exists
measurable set $A\subset (-1,1)$ of positive Lebesgue measure
such that $f>0$
a.e. on $A$
and, moreover,
$f\geq 0$ a.e. in $(-1,1)$.
Hence 
by positivity of Green's function
obtained in Lemma~\ref{lem:positivity}
and by Theorem~\ref{thm:conv}, we have
$$
    u(x) = \int_{-1}^1 G(x,y) f(y)\,\mathrm{d}y\geq
    \int_A G(x,y) f(y)\,\mathrm{d}y>0\,,
$$
for any $x\in (-1,1)$.
\hfill$\blacksquare$
\par\medskip

\section{Concluding remarks}
\label{sec:concrem}

\begin{table}[ht]
    \centering
    \begin{tabular}{|c|c|c|c|}
        \hline
        Result     & Assumptions on $f$ & $p$ & Result type \\
        \hline
        \hline
        Thm.~\ref{thm:C1:reg}  & {\rm(HF)} & $p > 1$ &  $C^1$-regularity. \\
        \hline
        Thm.~\ref{thm:boundness:h} & {\rm(HF)} & $p > 1$  & A~priori bound on $\|u\|_{\infty}$.\\
        \hline
        Thm.~\ref{thm:existence} &   $\displaystyle \|f\|_{L^1(-1,1)} < H (\sin\varphi)^{p - 1}$& 
        $p > 1$ & Existence of weak solution.\\
        \hline
        Lem.~\ref{prop:nonlinear:to:linear} & {\rm(HF)} & $p > 2$ & Linearization.\\
        \hline
        Thm.~\ref{thm:boundedness:derivative} & {\rm(HF)} & $p > 2$ & Bound on $\|u'\|_{\infty}$.\\
        \hline
        Thm.~\ref{thm:apriory:bound:derivative} & \makecell{{\rm(HF)}; \\ $\displaystyle \|f\|_{L^1(-1,1)} \leq \beta < H(\sin\varphi)^{p - 1}$ }& $p > 2$ & A~priori bound on $\|u'\|_{\infty}$.\\
        \hline
        Prop.~\ref{wmp} & $f\in L^1(-1,1)$; $f \geq 0$ & $p > 2$ & Weak Max. Principle.\\
        \hline
        Thm.~\ref{smp} & $f\in L^1(-1,1)$; $f \geq 0$; $f\not\equiv0$ & $p > 2$ & Strong Max. Principle.\\
        \hline
    \end{tabular}
    \caption{Main results.}
\label{tab:thm:overview}
\end{table}

We proposed a new mathematical model of groundwater flow over inclined impermeable bed. It turned out that this lead us to a strongly nonlinear problem, which is quite difficult to investigate. For summary of our main results, see Table~\ref{tab:thm:overview} together with their assumptions on $f$ and $p$.

The most important of our qualitative results cover the case $p>2$ only, which corresponds to flow in media of low permeability such as certain sandstones, fine sands, clays or certain types of soils (see, e.g., \cite{King1898, SoniIslamBasak1978, Zunker1920}). 
Let us note that the case $3/2<p<2$ is very important in applications too, as it corresponds to flows in coarse grained porous media such as gravels, see, e.g., \cite{BenediktGirgKotrlaTakac2018, SoniIslamBasak1978, Zunker1920}. Thus generalizations of results valid for $p>2$ stated in Table~\ref{tab:thm:overview} to include all $p>3/2$ would be of great importance. For further discussion, see Remark~\ref{rem:pl2}.

Let us also note that $f\in L^1(-1,1)$ and $f\geq 0$ imply (HF), so we do not explicitly assume (HF) in results concerning maximum principles, where $f\geq 0$ is natural assumption. We did not tackle Weak and Strong Comparison Principles, which for nonlinear problems are more difficult to prove than maximum principles. We suggest an investigation of comparison principles as a very promising direction
of research from both theoretical and practical point of view.

Another interesting direction of research would be to consider situations when (HF) is not satisfied. Then the solution $u$ can reach the value $-H$ in some subdomains of $(-1,1)$. This situation would correspond to the case that there is no groundwater over such subdomains.
This is, however, realistic scenario worth of further research. Last but not least, it is also of great importance to study the case when $H=0$, which corresponds to the case when the ditches are dry.

Finally, let us remark that the use of modern geophysical non-invasive techniques helps to understand internal structures of slopes, see, e.g., \cite{Duffek2023, Otto2006, SassWolny2001,  VolkelLeopoldRoberts2001}. Using these techniques, locations with porous media layered over impermeable beds (bedrock) have been observed in natural landscapes. These locations are encountered, e.g., in hilly or mountainous areas of central Europe, which shed water into the surrounding highly populated areas and are rich in precipitation. Hence, it is of great practical importance to understand water flow in such locations. But in these cases, boundary conditions other than Dirichlet need to be considered. 
Also the imaging techniques from \cite{Duffek2023,Otto2006,SassWolny2001} reveal that the bedrock does not have to be flat hyperplane and realistic model needs to take into account possibly curved surface of the bedrock. This suggests another interesting direction of future research.

\appendix

\section{}
\label{app:a}
\begin{comment}
To demonstrate that the operator $A$ defined by the left-hand side of~\eqref{eq:trunk} is bounded, coercive, and pseudomonotone, a convenient approach is to recognize that~\eqref{eq:trunk} is a specific instance of \cite[Eq.~(2.51), p.~44]{Roubicek2ndedit}.
\end{comment}

A convenient approach how to show that the operator $$A\colon W_0^{1,p}(-1,1)\to W^{-1,p'}(-1,1)$$ defined by the left-hand side of~\eqref{eq:trunk} has desired properties, is to recognize that \eqref{eq:trunk} is a specific instance of \cite[Eq.~(2.51), p.~44]{Roubicek2ndedit}. Then it suffices to write expressions for functions 
$a, c\colon (-1,1)\times
\mathbb{R}\times
\mathbb{R}$ from \cite[Eq.~(2.51), p.~44]{Roubicek2ndedit} to match left-hand side of~\eqref{eq:trunk}
and
verify structural assumptions
\cite[(2.54), (2.55), (2.65), and (2.92)]{Roubicek2ndedit}.

\begin{comment}   
A convenient approach how to show that the operator $A$ defined by~\eqref{A:def} is bounded, coercive, and pseudomonotone is to verify structural assumptions
\cite[(2.54), (2.55), (2.65), and (2.92)]{Roubicek2ndedit}
for $a, c, \colon (-1,1)\times
\mathbb{R}\times
\mathbb{R}$ defined in our case by
\end{comment}
Indeed, it is easy to see that the choice
\begin{eqnarray}
a(x,r,s)
&%\stackrel{\mathrm{def}}{=}
=& \left(H+T_{k}(r)\right)(\cos\varphi)^{p-1}
|s+\tan\varphi|^{p-2}(s+\tan\varphi)\,, \\
c(x,r,s)
&%\stackrel{\mathrm{def}}{=}
=& 0\,,
\end{eqnarray}
in \cite[Eq.~(2.51), p.~44]{Roubicek2ndedit}
matches the expression on the left-hand side of~\eqref{eq:trunk}.
Let us recall that
$T_k(t)=
\max\{-k, \min\{t, k\}\}$.
Since  
$\Gamma_{\mbox{N}}$ from
\cite[(2.54)]{Roubicek2ndedit}
is $\Gamma_{\mbox{N}}=\emptyset$ in
our case, the function $b$
from \cite[(2.54)]{Roubicek2ndedit} is set to be zero. Both functions $b=0$ and
$c=0$ trivially satisfy all required assumptions posed thereon in \cite{Roubicek2ndedit}. Thus
we do not mention these conditions on $b$ and $c$ explicitly here.
\par\medskip
\paragraph{\bf Boundedness.}
Taking into account that our function $a$ is continuous in variables $r$ and $s$ (and does not depend on $x$), it is a Carath\'{e}odory function.
From the growth conditions \cite[(2.55a--c)]{Roubicek2ndedit},
we need to verify only 
\cite[(2.55a)]{Roubicek2ndedit}, since the other are satisfied trivially in our situation. Taking account that
\begin{multline*}
|a(x,r,s)|=
\left|\vphantom{\int}
\left(H+T_{k}(r)\right)
(\cos\varphi)^{p-1}
|s+\tan\varphi|^{p-2}(s+\tan\varphi)
\right|
\leq \\
|H+\max\{-k, \min\{r, k\}\}|\,|s+\tan\varphi|^{p-1}
\leq 2^{p-2}(H+k) \left(|s|^{p-1}
+ |\tan\varphi|^{p-1}
\right)
\end{multline*}
is satisfied for all $x\in (-1,1)$, $r,s\in\mathbb{R}$, the condition
\cite[(2.55a)]{Roubicek2ndedit}
is verified.
According to \cite[Lemma~2.31]{Roubicek2ndedit},
conditions 
\cite[(2.54) and (2.55)]{Roubicek2ndedit} ensure that {\it the~operator $A\colon W_0^{1, p}(-1,1) \rightarrow W^{-1, p'}(-1,1)$ is bounded}.

\par\medskip
\paragraph{\bf Coercivity.}
Observe that
$$
\lim\limits_{s\to\pm\infty}
\frac{
|s+\tan\varphi|^{p-2}(s+\tan\varphi)\,s}{|s^p|}=1 
$$
for any $p>1$ and $\varphi\in (0,\pi/2)$.
Hence
$$
\lim\limits_{s\to\pm\infty}
|s+\tan\varphi|^{p-2}(s+\tan\varphi)\,s - \frac{1}{2}|s|^p
=+\infty
$$
and, by continuity argument,
$$
\mu
\stackrel{\mathrm{def}}{=}
-\min\left\{0, \min\limits_{s\in\mathbb{R}}|s+\tan\varphi|^{p-2}(s+\tan\varphi)\,s - \frac{1}{2}|s|^p\right\}\geq 0$$
exists and is finite.
Then
\begin{multline}
|s+\tan\varphi|^{p-2}(s+\tan\varphi)\,s = \\
\frac{1}{2}|s|^p +
|s+\tan\varphi|^{p-2}(s+\tan\varphi)\,s - \frac{1}{2}|s|^p \geq \\
\frac{1}{2}|s|^p +
\min\limits_{s\in\mathbb{R}}
\left(|s+\tan\varphi|^{p-2}(s+\tan\varphi)\,s - \frac{1}{2}|s|^p\right) \geq 
\frac{1}{2}|s|^p - \mu
\end{multline}
For $0<k<H$, $H+T_k(r)\geq H-k>0$. Thus
\begin{multline}
a(x,r,s)\, s = \left(H+T_{k}(r)\right)(\cos\varphi)^{p-1}
|s+\tan\varphi|^{p-2}(s+\tan\varphi)\,s 
\geq \\
\left(H-k\right)(\cos\varphi)^{p-1}
|s+\tan\varphi|^{p-2}(s+\tan\varphi)\,s
\geq \\
\frac{1}{2} \left(H-k\right)(\cos\varphi)^{p-1} |s|^p - \left(H-k)\right)(\cos\varphi)^{p-1} \mu\,,
\end{multline}
which means that the structural condition
\cite[(2.92a)]{Roubicek2ndedit}
is verified.
According to \cite[Lem.~2.35]{Roubicek2ndedit},
{\it the~operator $A\colon W_0^{1, p}(-1,1) \rightarrow W^{-1, p'}(-1,1)$ is coercive}.

\par\medskip
\paragraph{\bf Pseudomonotonicity.}
In order to prove pseudomonotonicity, we will make use of the following well-known results.
\begin{comment}
{\color{red} \verb|https://matematicas.uam.es/~ireneo.peral/ICTP.pdf| Lemma A.05, str. 80 + je tam odkaz na francouzsky clanek z roku 1979, kde by to melo byt. Poslu emailem.}
\end{comment}

\begin{proposition}{\rm(see, e.g., 
\cite[pp.~210--211]{Simon78}.)}
Let $p>1$ and $N\in\mathbb{N}$. Then
there exists a constant
$c_p>0$ (depending on $p$) such that,
for all $x, y \in \mathbb{R}^N$, the following inequality holds
$$
\left\langle|x|^{p-2} x-|y|^{p-2} y, x-y\right\rangle_{\mathbb{R}^N} \geq\left\{\begin{array}{cl}
c_p|x-y|^p & \text { if } p \geq 2\,, \\
c_p \frac{|x-y|^2}{(|x|+|y|)^{2-p}} & \text { if } 1<p<2\,.
\end{array}\right.
$$
\end{proposition}

Now we are ready to verify condition \cite[Eq.~(2.65)]{Roubicek2ndedit}
(the so called monotonicity in the main part).
For this we introduce new variables $\sigma=s-\tan\varphi$, $\tilde\sigma=
\tilde s-\tan\varphi$ and observe $s-\tilde s=\sigma-\tilde\sigma$. With this, we have
\begin{multline*}
(a(x, r, s)-a(x, r, \tilde{s})) \cdot(s-\tilde{s})= \\
(a(x, r, \sigma + \tan\varphi)-a(x, r, \tilde\sigma +\tan\varphi)) \cdot(\sigma-\tilde{\sigma})= \\
\left(H+T_{k}(r)\right)
(\cos\varphi)^{p-1}\left(\vphantom{\int}
|\sigma|^{p-2}(\sigma)
-
|\tilde{\sigma}|^{p-2}(\tilde{\sigma})
\right)\cdot(\sigma-\tilde{\sigma})
\geq \\ 
\left\{\begin{array}{cl}
c_p\left(H+T_{k}(r)\right)(\cos\varphi)^{p-1}|\sigma-\tilde{\sigma}|^p & \text { if } p \geq 2 \\[0.3cm]
\displaystyle
c_p\left(H+T_{k}(r)\right)
(\cos\varphi)^{p-1}\frac{|\sigma-\tilde{\sigma}|^2}{(|\sigma|+|\tilde{\sigma}|)^{2-p}} & \text { if } 1<p<2
\end{array}\right\}\geq 0\,,
\end{multline*}
for all $x\in(-1,1)$,
$r, s, \tilde{s}\in\mathbb{R}$.
Now, by \cite[Lem.~2.32]{Roubicek2ndedit},
conditions 
\cite[(2.54), (2.55), and (2.65)]{Roubicek2ndedit} ensure that {\it the~operator $A\colon W_0^{1, p}(-1,1) \rightarrow W^{-1, p'}(-1,1)$ is pseudomonotone}.

\section*{In Memoriam: Prof. Neuberger's Legacy of Inspiration and Guidance}
As a young researcher attending the ``Variational Methods: Open Problems, Recent Progress, and Numerical Algorithms" conference held in Flagstaff AZ, USA, in 2002, Petr Girg had the privilege of meeting Professor John~W. Neuberger and engaging in several stimulating discussions that brought to his attention a deeper connection between theoretical, applied, and numerical mathematics. These discussions had a significant impact on Petr's later career, and he would like to take this opportunity to express his sincere gratitude for Professor Neuberger's insightful advice and encouragement. Professor Neuberger's absence will be deeply felt by younger researchers across the field of mathematics.

\section*{Acknowledgements}
P.~Girg and L. Kotrla were supported by the~Grant Agency of~the~Czech
Republic, Grant No.~22-18261S.

\bibliographystyle{acm}
\bibliography{ref.bib}

\end{document}